\documentclass[12pt,oneside]{amsart}
\usepackage{cite}
\usepackage{latexsym}
\usepackage{layout}
\usepackage{amsfonts}
\usepackage{amssymb}
\usepackage{epsfig}

\newtheorem{thm}{Theorem}[section]

\newtheorem{lem}[thm]{Lemma}

\newtheorem{prp}[thm]{Proposition}

\begin{document}
\title{Quadratic Residues and Non--residues in Arithmetic Progression}
\author{Steve 
wright} 
\address{Steve Wright, Department of Mathematics and Statistics,
Oakland University, Rochester, MI 48309. Email: wright@oakland.edu}

\begin{abstract}
Let $\mathcal{S}$ be an infinite set of non-empty, finite subsets of the nonnegative integers. If $p$ is an odd prime, let $c(p)$ denote the cardinality of the set $ \{ S \in \mathcal{S} : S \subseteq \{1,\dots,p-1\}$ and $S$ is a set of quadratic residues (respectively, non-residues) of $p$\}. 
When $\mathcal{S}$ is constructed in various ways from the set of all arithmetic progressions of nonnegative integers, we determine the sharp asymptotic behavior of $c(p)$ as $p \rightarrow +\infty$. Generalizations and variations of this are also established, and some problems connected with these results that are worthy of further study are discussed. 
 \end{abstract}
\maketitle
\markboth{}{}
\noindent \textit{keywords}: \textrm{quadratic residue, quadratic non-residue, arithmetic progression, asymptotic approximation, Weil sum}

\section{Introduction}
\label{intro}
If $p$ is an odd prime, an integer $z$ is said to be a \emph{quadratic residue} (respectively, \emph{quadratic non-residue}) of $p$ if the equation $x^2 \equiv z $ mod $ p$ has (respectively, does not have) a solution $x$ in integers. It is a theorem going all the way back to Euler that exactly half of the integers from 1 through $p-1$ are quadratic residues of $p$, and it is a fascinating problem to investigate the various ways in which these residues are distributed among 1, 2,\dots, $p-1$. In this paper, our particular interest lies in studying the problem of the distribution of residues and non-residues among the arithmetic progressions which can occur in the set $\{1, 2, \dots, p-1\}$.

We begin with a litany of notation and terminology that will be used systematically throughout the rest of this paper. If $m \leq n$ are integers, then $[m,n]$ will denote the set of all integers that are at least $m$ and no greater than $n$, and $[m,+\infty)$ will denote the set of all integers that exceed $m-1$. For any odd prime $p$, we let $R(p)$ (respectively, $NR(p)$) denote the set of all quadratic residues (respectively, non-residues) of $p$ in the interval $[1,p-1]$. If $\{a(p)\}$ and $\{b(p)\}$ are sequences of real numbers defined for all primes $p$ in an infinite set $S$, then we will say that $a(p)$ is (sharply) \emph{asymptotic to} $b(p)$ \emph{as}  $p \rightarrow +\infty$ \emph{inside} $S$, denoted as $a(p) \sim b(p)$, if
\[
\lim_{\substack{p\rightarrow +\infty \\ p\in S}} \frac{a(p)}{b(p)} = 1,
\]
and if $S = [0,+\infty)$, we simply delete the phrase "inside $S$". If $x$ is a real number, $[x]$ will denote the greatest integer that does not exceed $x$. Finally, if $A$ is a set then $|A|$ will denote the cardinality of $A$, $ 2^A$ will denote the set of all subsets of $A$, $\mathcal{E}(A)$ will denote the set of all nonempty subsets of $A$ of even cardinality, and $\emptyset$ will denote the empty set. We also note once and for all that $p$ will always denote a generic odd prime.

 Our work here, in both spirit and method, has its origins in some classical results of H. Davenport. In the papers [5, 6, 7], Davenport considers the problem of estimating the number $R_s(p)$ (respectively, $N_s(p)$) of sets of $s$ consecutive quadratic residues (respectively, non-residues) of an odd prime $p$ that occur inside $[1,p-1]$. The expected number is about $2^{-s}p$, and in [7, Corollary of Theorem 5], Davenport showed that as $p\rightarrow+\infty$, both $R_s(p)$ and $N_s(p)$ are asymptotic to $2^{-s}p$. His method of proof is based on the following clever idea. If $Z_p$ is the field of $p$ elements, then the Legendre symbol of $p$ defines a real (primitive) multiplicative character $\chi_p : Z_p \rightarrow [-1,1]$ on $Z_p$. If $\varepsilon \in \{-1,1\}$, form the sum
 \begin{equation*}
 2^{-s}\sum_{x=1}^{p-s}\  \prod_{i=0}^{s-1} \Big(1+\varepsilon \chi_p(x+i)\Big)\tag{1.1}
 \end{equation*} 
 and note that the value of this sum is $R_s(p)$ (respectively, $N_s(p)$) when $\varepsilon = 1$ (respectively, $\varepsilon = -1$). Davenport rewrote this sum as
\begin{equation*}
 2^{-s}(p-s) +2^{-s} \sum_{\emptyset\ \not  = \ T\ \subseteq\ [0,s-1]} \varepsilon^{|T|}\Big (\sum_{x=1}^{p-s} \chi_p\Big(\prod_{i \in T} (x+i)\Big)\Big)\tag{1.2} 
\end{equation*}
and then used the theory of Hasse $L$-functions to prove that there exists positive absolute constants $\sigma <1$ and $C$ such that for all $p$ sufficiently large,
\begin{equation*}
\left|\sum_{x=1}^{p-s} \chi_p\Big(\prod_{i\in T} (x+i)\Big)\right|\leq Csp^\sigma .\tag{1.3}
\end{equation*}
When this estimate is applied in (1.2), it follows immediately that $R_s(p)$  and $N_s(p)$ are asymptotic to $2^{-s}p$, with an asymptotic error which does not exceed $Csp^\sigma$, for all $p$ sufficiently large. In 1945, as a consequence of his landmark work on arithmetic algebraic geometry [15], A. Weil showed that the estimate (1.3) held with $C=2$ and $\sigma=1/2$, which produces an essentially optimal estimate for the asymptotic error. More generally, if $\chi_p$ is replaced in the sum in (1.3) by an arbitrary non-principal multiplicative character $\chi$ on a finite field $F$, Weil's work implies that the resulting sum also satisfies this improved estimate. Consequently, if $f$ is a polynomial over $F$, then sums of the form
\[
\sum_{x\in F} \chi(f(x))
\] 
 are called \emph{Weil sums}, a terminology to which we will adhere in this paper. For further applications, refinements, and discussion of Davenport's method, we refer the reader to [2, 3, 8] and especially [9, chapter 9].
 
 Our point of departure from Davenport's work is to notice that the sequence $\{x,x+1,\dots, x+s-1\}$ of $s$ consecutive positive integers that appears in (1.1)  is an instance of the sequence$\{x,x+b,\dots,x+b(s-1)\}$, an arithmetic progression of length $s$ and common difference $b$, with $b=1$. Thus, if $(b,s)\in [1,+\infty)\times[1,+\infty)$, and we set
 \[
 AP(b;s)=\Big\{ \{n+ib: i\in [0,s-1]\}: n\in [0,+\infty)\Big\},
 \]
the family of all arithmetic progressions of length $s$ and common difference $b$, it is natural to inquire about the asymptotics as $p\rightarrow +\infty$ of the number of elements of $AP(b;s)$ that are sets of quadratic residues (respectively, non-residues) of $p$ that occur inside $[1,p-1]$. We also consider the following related question: if $a\in [0,+\infty)$, set
\[
AP(a,b;s)=\Big\{ \{a+b(n+i):i\in [0,s-1]\}: n\in [0,+\infty)\Big\},
\]
the family of all arithmetic progressions of length $s$ taken from a fixed arithmetic progression 
\[
AP(a, b)= \{a+bn: n\in [0,+\infty)\}.
\]
We then ask for the asymptotics of the number of elements of $AP(a,b;s)$ that  are sets of quadratic residues (respectively, non-residues) of $p$ that occur inside $[1,p-1]$. 
Solutions of these problems will provide interesting insights into how often quadratic residues and non-residues appear as arbitrarily long arithmetic progressions.

 We will in fact consider the following generalization of these questions. For each $m\in [1,+\infty)$, let
 \[
 \textbf{a}=(a_1,\dots,a_m)\  \textrm{and}\ \textbf{b}=(b_1,\dots,b_m)
 \] 
 be $m$-tuples of nonnegative integers such that $(a_i,b_i)\not=(a_j,b_j)$, for all $i\not=j$. When the $b_i$'$s$ are distinct and positive, we set
 \[
 AP(\textbf{b}; s)=\Big\{\bigcup_{j=1}^{m} \{n+ib_j:i\in [0,s-1]\}:n\in [1, +\infty)\Big\},
 \]
 and when the $b_i$'s are all positive, we set
 \[
 AP(\textbf{a},\textbf{b};s)=\Big\{\bigcup_{j=1}^{m} \{a_j+b_j(n+i):i\in [0,s-1]\}:n\in [0,+\infty)\Big\}.
 \] 
 If $m=1$ then we recover our original sets $AP(b;s)$ and $AP(a,b;s)$. We now pose
 \begin{quote}
 \textbf{Problem 1} (respectively, \textbf{Problem 2}): determine the asymptotics as $p\rightarrow+\infty$ of the number of elements of $AP(\textbf{b};s)$ (respectively, $AP(\textbf{a}, \textbf{b};s)$) that are sets of quadratic residues of $p$ inside $[1,p-1]$.
 \end{quote}
 We also pose as \textbf{Problem 3} and \textbf{Problem 4} the problems which result when the phrase "quadratic residues" in the statements of Problems 1 and 2 is replaced by the phrase "quadratic non-residues".
 
In section 2 of this paper, we solve Problems 1 and 3. In fact, we will achieve considerably more than that; see our first main result, Theorem 2.3, and the definitions which precede it in section 2 for a precise statement of what we establish. A special case of Theorem 2.3, the one which gives a specific solution to Problems 1 and 3, asserts that if $\gamma$ denotes the cardinality of the set
\[
\bigcup_{j=1}^{m}\ \{ib_j: i\in [0, s-1]\},
\]
then the number of elements of the set $AP(\textbf{b}; s)\cap 2^{[1, p-1]}$ that are sets of quadratic residues (respectively, non-residues) of $p$ is asymptotic to $2^{-\gamma}p$ as $p\rightarrow +\infty$. 

Our approach to the proof of Theorem 2.3 tailors Davenport's method to this situation, i.e., we consider an appropriate sum of products, in analogy to the sum (1.1), whose value gives the number of elements of $AP(\textbf{b};s)$ that are sets of quadratic residues (respectively, non-residues) of $p$ inside $[1,p-1]$. We then expand the products in the sum and interchange the order of summation in what results to obtain an analog of the sum (1.2). This sum will consist of a dominant term that is a non-constant linear function of $p$ plus a remainder. In exact analogy with Davenport's argument, an estimate of the remainder term on the order of $O(\sqrt{p})$ will allow the dominant term to determine the asymptotics of the desired sequence. The reminder will have summands that are Weil sums for $\chi_p$ that are similar to the Weil sums which appear in (1.2), and can hence be estimated in the required manner by standard techniques.

In sections 3-7, we attack Problems 2 and 4 for $AP(\textbf{a},\textbf{b};s)$. Because of certain arithmetical interactions which can take place between the elements of the sets in $AP(\textbf{a},\textbf{b};s)$, the asymptotic behavior as $p\rightarrow +\infty$ of the number of elements of $AP(\textbf{a},\textbf{b};s)\cap 2^{[1, p-1]}$ which are sets of quadratic residues (respectively, non-residues) of $p$ is somewhat more complicated than what occurs for $AP(\textbf{b}; s)$. 

In order to explain the situation, we set
\begin{equation*}
q_\varepsilon(p)=|\{A\in AP(\textbf{a}, \textbf{b}; s)\cap 2^{[1, p-1]}: \chi_p(a)=\varepsilon, \ \textrm{for all}\ a\in A\}|
\end{equation*} 
and note that the value of $q_\varepsilon(p)$ for $\varepsilon=1$ (respectively, $\varepsilon=-1$) counts the number of elements of $AP(\textbf{a}, \textbf{b}; s)$ that are sets of quadratic residues (respectively, non-residues) of $p$ that are located inside $[1, p-1]$. As we mentioned before, it  will transpire that the asymptotic behavior of $q_\varepsilon(p)$ depends on certain arithmetic interactions that can take place between the elements of $AP(\textbf{a}, \textbf{b}; s)$. In order to see how this goes, first consider the set $B$ of \emph{distinct} values of the coordinates $b_1,\dots, b_m$ of $\textbf{b}$. If we declare the coordinate $a_i$ of $\textbf{a}$ and the coordinate $b_i$ of $\textbf{b}$ to \emph{correspond} to each other, then for each $b\in B$, we let $A(b)$ denote the set of all coordinates of $\textbf{a}$ which correspond to $b$. We then relabel the elements of $B$ as $b_1,\dots,b_k$, say, and for each $i\in [1, k]$, set $A_i=A(b_i)$,  \[
 S_i=\bigcup_{a\in A_i} \{a+b_il: l\in [0, s-1]\}, 
 \]
 and then let
 \[
 \alpha=\sum_i\ |S_i|,\   b=\max \{b_1,\dots, b_k\}.
 \]
 
 Next, suppose that
 \begin{quote}
 $(*)$  if $(i, j)\in [1, k]\times [1, k]$ with $i\not= j$ and $(a, a')\in A_i\times A_j$, then either $b_ib_j$ does not divide $a'b_i-ab_j$ or  $b_ib_j$ divides $a'b_i-ab_j$ with a quotient that exceeds $s-1$ in modulus.
 \end{quote}
We will then prove that as $p\rightarrow +\infty,\  q_\varepsilon(p)$ is asymptotic to $(b\cdot 2^\alpha)^{-1}p$. On the other hand, if the assumption $(*)$ does not hold then we show that the asymptotic behavior of $q_\varepsilon(p)$ falls into two distinct regimes, with each regime determined in a certain manner by the integral quotients
\begin{equation*}
\frac{a'b_i-ab_j}{b_ib_j},\  (a, a')\in  A_i\times A_j,\tag{**}
\]
whose moduli do not exceed $s-1$. More precisely, these quotients determine a positive integer $e<\alpha$ and a collection $\mathcal{S}$ of nonempty subsets of $[1, k]$ such that each element of  $\mathcal{S}$ has even cardinality and for which the following two alternatives hold:

$(i)$ if $\prod_{i\in S} b_i$ is a square for all $S\in \mathcal{S}$, then as $p\rightarrow +\infty,\  q_\varepsilon(p)$ is asymptotic to $(b\cdot 2^{\alpha-e})^{-1}p$, or

$(ii)$ if there is an $S\in \mathcal{S}$ such that $\prod_{i\in S} b_i$ is not a square, then there exist two disjoint, infinite sets of primes $\Pi_+$ and $\Pi_-$  whose union contains all but finitely many of the primes and such that $ q_\varepsilon(p)=0$ for all $p\in \Pi_-$, while as $p\rightarrow +\infty$ inside $\Pi_+$, $q_\varepsilon(p)$ is asymptotic to $(b\cdot 2^{\alpha-e})^{-1}p$.

Thus we see that when $(*)$ does not hold and $p\rightarrow +\infty$, either $ q_\varepsilon(p)$ is asymptotic to $(b\cdot 2^{\alpha-e})^{-1}p$ or $q_\varepsilon(p)$ asymptotically oscillates infinitely often between 0 and $(b\cdot 2^{\alpha-e})^{-1}p$.

In light of what we have just discussed, it is no surprise that the solution of Problems 2 and 4 for  $AP(\textbf{a},\textbf{b};s)$ will involve  a bit more effort than the solution of Problems 1 and 3 for  $AP(\textbf{b}; s)$. In order to analyze the asymptotic behavior of $q_\varepsilon(p)$, we follow the same strategy as before: using an appropriate sum of products involving $\chi_p,\  q_\varepsilon(p)$ is expressed as a sum of a dominant term and a remainder. If the dominant term is a non-constant linear function of $p$ and the remainder term is $O(\sqrt p\log p)$, then the asymptotic behavior of $ q_\varepsilon(p)$ will be in hand.    
 
We in fact will implement this strategy when the set  $AP(\textbf{a}, \textbf{b}; s)$ in the definition of $q_\varepsilon(p)$ is replaced by a slightly more general set; for a precise statement of what we establish, see Theorem 6.1 in section 6, the second principal result of this article. We then deduce our results for $AP(\textbf{a}, \textbf{b}; s)$ from this more general result in section 7, where, in particular, the reader can find the precise manner in which the integral quotients $(**)$ whose moduli do not exceed $s-1$ determine the parameter $e$ and collection of sets $\mathcal{S}$ discussed above. Section 3 contains the required estimate of the remainder term and a preliminary calculation of the required dominant term in this more general situation. The dominant term that arises here is considerably more complicated than the one which occurs for $AP( \textbf{b}; s)$, and so as an auxiliary to our analysis of it, we define and study a device in section 4, the $(B, \textbf{S})$-\emph{signature} of a prime, which is then used in section 5 to finish the calculation of the dominant term. In section 8, we study an interesting class of $2k$-tuples $(\textbf{a}, \textbf{b})$ for which the parameters $\alpha$ and $e$ can be easily calculated, and we use this fact to illustrate how Theorem 6.1 determines the asymptotic behavior of $q_\varepsilon(p)$ in concrete situations. Our discussion concludes in section 9 with remarks on some problems worthy of further study which arise naturally from the work done in this paper. 

\section{the results for $AP(\textbf{b}; s)$}
We begin this section with some terminology and notation that will allow us to state our results precisely and concisely. Let $Z=\{z_1,\dots,z_r\}$ be a finite subset of $[0,+\infty)$ with its elements indexed in increasing order $z_i<z_j$ for $ i<j$. 
We let
\[
\mathcal{S}(Z)=\Big\{\{n+z_i:i\in [1,r]\}:n\in [1,+\infty)\Big\},
\]
the set of all shifts of $Z$ to the right by a positive integer. Let $\varepsilon$ be a choice of signs for $[1,r]$, i.e., a function from $[1,r]$ into $\{-1,1\}$. If $S=\{n+z_i:i\in [1,r]\}$ is an element of $\mathcal{S}(Z)$, we will say that the triple $(S,\varepsilon,p)$ is a \emph{residue pattern of p} if
\[
\chi_p(n+z_i)=\varepsilon(i),\  \forall\  i\in [1,r].
\]
The set $\mathcal{S}(Z)$ has the \emph{universal pattern property} if there exists $p_0>0$ such that for all $p\geq p_0$ and for all choices of signs $\varepsilon$ for $[1,r]$, there is a set $S\in \mathcal{S}(Z)\cap 2^{[1,p-1]}$ such that $(S,\varepsilon,p)$ is a residue pattern of $p$. $\mathcal{S}(Z)$ hence has the universal pattern property if and only if for all $p$ sufficiently large, $\mathcal{S}(Z)$ contains a set that exhibits any fixed but arbitrary pattern of quadratic residues and non-residues of $p$. This property is inspired directly by Davenport's work: using this terminology, we can state the result of [7, Corollary of Theorem 5] for quadratic residues as asserting that if $s\in [1,+\infty)$ then $\mathcal{S}([0,s-1])$ has the universal pattern property, and moreover, for any choice of signs $\varepsilon$ for $[1,s]$, the cardinality of the set
\[
\{S\in \mathcal{S}([0,s-1])\cap 2^{[1,p-1]}: (S,\varepsilon,p)\  \textrm{is a residue pattern of}\ p\}
\]    
is asymptotic to $2^{-s}p$ as $p\rightarrow +\infty$. Note that if $\varepsilon$ is the choice of signs that is either identically 1 or identically $-1$ on $[1,s]$, then we recover the results that were discussed in section 1.

Suppose now that there exists nontrivial gaps between elements of $Z$, i.e., $z_{i+1}-z_i\geq 2$ for at least one $i\in [1,r-1]$. It is then natural to search for elements $S$ of $\mathcal{S}(Z)$ such that the quadratic residues (respectively, non-residues) of $p$ inside $[\min{S},\max{S}]$ consists precisely of the elements of $S$, so that $S$ acts as the "support" of quadratic residues or non-residues of $p$ inside the minimal interval of consecutive integers containing $S$. We formalize this idea by declaring $S$ to be a \emph{residue} (respectively,  \emph{non-residue}) \emph{support set of p} if $S=R(p)\cap [\min S,\max S]$ (respectively, $S=NR(p)\cap [\min S,\max S]$). We then define $\mathcal{S}(Z)$ to have the \emph{residue} (respectively, \emph{non-residue}) \emph{support property} if there exist $p_0>0$ such that for all $p\geq p_0$, there is a set $S\in \mathcal{S}(Z)\cap 2^{[1,p-1]}$ such that $S$ is a residue (respectively, non-residue) support set of $p$. 

We now use Davenport's method to establish the following proposition, which generalizes [7, Corollary of Theorem 5] for quadratic residues.

\begin{prp}
\label{prp1}
If $Z$ is any nonempty, finite subset of $[0,+\infty)$, then $\mathcal{S}(Z)$ has the universal pattern property and both the residue and non-residue support properties. Moreover, if $\varepsilon$ is a choice of signs for $[1,|Z|]$,
\[
c_\varepsilon(Z)(p)=\left| \{S\in \mathcal{S}(Z)\cap 2^{[1,p-1]}: (S,\varepsilon,p)\textrm{is a residue pattern of}\  p\}\right|, and
\]
\[
c_\sigma(Z)(p)=\left| \{S\in \mathcal{S}(Z)\cap 2^{[1,p-1]}: S\textrm{is a residue (respectively, non-residue) support set of}\  p\}\right|,
\]
then as $p\rightarrow +\infty$,
\[
 c_\varepsilon(Z)(p)\sim 2^{-|Z|}p\ and\ c_\sigma(Z)(p)\sim 2^{-(1+\max Z-\min Z)}p.
 \]
 \end{prp} 
 
\textit{Proof}. Suppose that the asserted asymptotics of $c_\varepsilon(Z)(p)$ has been established for all nonempty, finite subsets $Z$ of $[0,+\infty)$. Then the asserted asymptotics for $c_\sigma(Z)(p)$ can be deduced from that by means of the following trick. Let $Z\subseteq [0,+\infty)$ be nonempty and finite. Define the choice of signs $\varepsilon$ for [min $Z$, max $Z$] to be 1 on $Z$ and $-1$ on [min $Z$, max $Z$]$\setminus$$Z$. Now for each $p$, let
\[
\mathcal{S}(p)=\{S\in \mathcal{S}(Z)\cap 2^{[1,p-1]}: S\ \textrm{is a residue support set of}\  p\},
\]
\[ 
\mathcal{R}(p)=\{S\in \mathcal{S}([\min Z,\max Z])\cap 2^{[1,p-1]}: (S,\varepsilon,p)\textrm{is a residue pattern of}\  p\}.
\]    
If to each $E\in \mathcal{R}(p)$ (respectively, $F\in \mathcal{S}(p)$), we assign the set $f(E)=E\cap R(p)$ (respectively, $g(F)=[\min F,\max F]$), then $f$ (respectively, $g$) maps $\mathcal{R}(p)$ (respectively, $\mathcal{S}(p)$) injectively into $\mathcal{S}(p)$ (respectively, $\mathcal{R}(p)$). Hence $\mathcal{R}(p)$ and $\mathcal{S}(p)$ have the same cardinality.
Because of our assumption concerning the asymptotics of $c_\varepsilon([\min Z,\max Z])(p)$, it follows that as $p\rightarrow +\infty$,
\[
c_\sigma(Z)(p)=\left| \mathcal{S}(p)\right|=\left| \mathcal{R}(p)\right|\sim 2^{-\left| [\min Z,\max Z]\right|}\ p=2^{-(1+\max Z-\min Z)}\ p.
\]
This establishes the conclusion of the proposition with regard to residue support sets, and the conclusion with regard to non-residue support sets follows by repeating the same reasoning after $\varepsilon$ is replaced by $-\varepsilon$.

If $\varepsilon$ is now an arbitrary choice of signs for$[1,|Z|]$, it hence suffices to deduce the asserted asymptotics of $c_\varepsilon(Z)(p)$. Letting $r(p)=p-\max Z-1$, we have for all $p$ sufficiently large that
\[
c_\varepsilon(Z)(p)=2^{-|Z|}\sum_{x=1}^{r(p)}\ \prod_{i=1}^{|Z|} \Big(1+\varepsilon(i)\chi_p(x+z_i)\Big).
\]
This sum can hence be rewritten as
\[
2^{-|Z|}(p-\max Z)+2^{-|Z|}\sum_{\emptyset\ \neq\ T\ \subseteq\ [1,|Z|]}\  \prod_{i\in T} \varepsilon(i)\Big(\sum_{x=1}^{r(p)} \chi_p\Big(\prod_{i\in T} (x+z_i)\Big).
\] 
The asserted asymptotics for $c_\varepsilon(Z)(p)$ now follows from an application of the next lemma to the Weil sums in the second term of this expression. $\Box$                                 

The following lemma provides the estimate of Weil sums that will be required in this section and the next. It is valid in much greater generality, but we state it here in a form that is most convenient for the work done in this paper. A proof can be had by combining the completion method for estimation of character sums as set forth in [11, section 12.2] with well-known estimates of hybrid or mixed Weil sums that are available from [14, Theorem 2.2G] or [12, 13]. 
\begin{lem}
\label{lem2}
There exists $M>0$ such that the following statement is true: if $p\geq M$, if $f\in Z_p[x]$ is a monic polynomial over $Z_p$ of degree $d\geq 1$ with distinct roots in $Z_p$, and if $N\in [0,p-1]$, then
\begin{equation*}
\left| \sum_{x=0}^{N} \chi_p(f(x))\right|\leq 2d\sqrt p\log p.
\end{equation*}
\end{lem}

The following theorem is one of the principal results of this paper. In particular, if the choice of signs $\varepsilon$ in part ($i$) of the theorem is taken to be either identically 1 or identically $-1$, we obtain the solution of Problems 1 and 3 that were posed for $AP(\textbf{b}; s)$ in section 1.
\begin{thm}
\label{thm3}
 If $\textbf{b}=(b_1,\dots,b_k)\in [1,+\infty)^k$ with $b_i\not= b_j$ for $i\not= j$, then $AP(\textbf{b}; s)$ has the universal pattern property and both the residue and non-residue support properties. Moreover, if $b=\max_i\{b_i\}$,
\[
\gamma=\left|\bigcup_{j=1}^{k}\ \{ib_j: i\in [0, s-1]\}\right|,
\]
$\varepsilon$ is a choice of signs for $[1, \alpha]$,
\[
c_\varepsilon(p)=|\{S\in AP(\textbf{b}; s)\cap 2^{[1, p-1]}: (S,\varepsilon, p)\textit{ is a residue pattern of p}\}|, and
\]
\[
c_\sigma(p)=|\{S\in AP(\textbf{b}; s)\cap 2^{[1, p-1]}: S\ \textit{is a residue (respectively, non-residue) support set of p}\}|,
\]
then as $p\rightarrow+\infty$,
\[
c_\varepsilon(p)\sim 2^{-\gamma}p\  and\ c_\sigma(p)\sim 2^{-(1+b(s-1))}p.
\]

\end{thm}

\emph{Proof}. If $Z$ is the union whose cardinality we have set equal to $\gamma$ in the statement of the theorem, then
\[
AP(\textbf{b}; s)=\mathcal{S}(Z),
\]
and so the conclusion of this theorem is an immediate consequence of Proposition 2.1.  $\  \Box$

\emph{Remarks}. (1) A careful inspection of the proof of Theorem 2.3 reveals that one can replace in that proof the incomplete Weil-sum estimates of Lemma 2.2 by the complete Weil-sum estimate
\[
\left| \sum_{x=0}^{p-1} \chi_p(f(x))\right|\leq 2d\sqrt p,
\]
where $f$ is any polynomial as specified in the statement of Lemma 2.2 and $d$ is the degree of  $f$  ([15], [9, section 9.4]). It is then straightforward to deduce from the proof of Theorem 2.3 the following error estimates for the asymptotic approximations in that theorem:  for all $p$ sufficiently large,
\[
|c_\varepsilon(p)-2^{-\gamma}p|\leq 2\gamma\sqrt p,
\]
\[
|c_\sigma(p)-2^{-(1+b(s-1)}p|\leq 2(1+b(s-1))\sqrt p.
\]

(2)  Suppose that $k\geq 2$, $b_i<b_j$ for $i<j$, and the greatest common divisor of $b_i$ and $b_j$ is 1 for $i\not= j$. In this situation, there is an elegant formula for the parameter $\gamma$ in the statement of Theorem 2.3, which we will now derive. Let $A_j=\{ib_j: i\in [0, s-1]\}$ and then use the principle of inclusion and exclusion to conclude that
\begin{equation*}
\gamma=\Big|\bigcup_{j} A_j\Big|=\sum_{l=1}^{k} (-1)^{l+1}\left(\sum_{T\subseteq [1, k],\ |T|=l} \Big|\bigcap_{i\in T} A_i\Big|\right).\tag{2.1}
\end{equation*}
Because the $b_i$'$s$ are pairwise relatively prime, the theory of linear Diophantine equations implies that
\[
\bigcap_{i\in T} A_i=\{m\prod_{i\in T} b_i: m\in [0, \left[\frac{s-1}{\max_{i\in T}\Big\{\prod_{j\in T\setminus \{i\}} b_j\Big\}}\right]]\}.
\] 
Since $b_i<b_j$ for $i<j$,
\[
\max_{i\in T}\left\{\prod_{j\in T\setminus \{i\}} b_j\right\}=\prod_{i\in T\setminus \{\min T\}}b_i,
\]
where we assign the value 1 to any empty product occurring here, and so
\[
\Big|\bigcap_{i\in T} A_i\Big|=1+\left[\frac{s-1}{\prod_{i\in T\setminus \{\min T\}}b_i}\right].
\]   
Hence from (2.1) it follows that
\begin{equation*}
\gamma=1+k(s-1)+\sum_{l=2}^{k} (-1)^{l+1}\Big(\sum_{T\subseteq [1, k],\ |T|=l}\left[\frac{s-1}{\prod_{i\in T\setminus \{\min T\}}b_i}\right]\Big).\tag{2.2}
\end{equation*}

We want to pin down a bit more precisely the sum in the parentheses on the right-hand side of (2.2). Fix $l\in [2, k]$, let $\mathcal{S}_i=\{T\subseteq [1, k]: |T|=i\}$
and consider the map $\varphi: \mathcal{S}_l\rightarrow \mathcal{S}_{l-1}$ defined by $\varphi(T)=T\setminus \{\min T\}$. Then $\varphi(\mathcal{S}_l)=\{S\subseteq [2, k]: |S|=l-1\}$, and if $S\in \varphi(\mathcal{S}_l)$ then $\varphi^{-1}(S)= \{\{x\}\cup S: x\in [1, \min S-1]\}$, hence $ |\varphi^{-1}(S)|=\min S-1$. Because $\mathcal{S}_l$ is the pairwise disjoint union
\[
\bigcup_{S\subseteq [2, k],\ |S|=l-1} \varphi^{-1}(S),
\]
it follows that
\begin{equation*}
\sum_{T\in \mathcal{S}_l} \left[\frac{s-1}{\prod_{i\in T\setminus\{\min T\}}b_i}\right]=\sum_{S\subseteq [2, k],\ |S|=l-1} (\min S-1)\left[\frac{s-1}{\prod_{i\in S}b_i}\right].\tag{2.3}
\end{equation*}
If we now define  $R(l, m)$ to be the set 
\[
\{S\subseteq [2, k]: |S|=l-1, \min S=m\},\ (l, m)\in [2, k]\times[2, k-l+2],
\]
then $\{S\subseteq [2, k]: |S|=l-1\}$ is the pairwise disjoint union 
\[
\bigcup_{m\in [2, k-l+2]} R(l, m),
\] 
and so it follows that
\begin{equation*}
\textrm{sum on the left-hand side of (2.3)}=\sum_{m=2}^{k-l+2} (m-1)\sum_{T\in R(l, m)}\left[\frac{s-1}{\prod_{i\in T} b_i}\right].\tag{2.4}
\end{equation*}
We conclude from (2.2), (2.3), and (2.4) that 
\[
\gamma=1+k(s-1)+\sum_{l=2}^{k}\ \ (-1)^{l+1}\sum_{m=2}^{k-l+2} (m-1)\sum_{T\in R(l, m)} \left[\frac{s-1}{\prod_{i\in T} b_i}\right].
\]
 
\section{The results for $AP(\textbf{a}, \textbf{b}; s)$: preliminaries} 
 Let  $(m, s)\in [1, +\infty)\times [1, +\infty)$ and let $(\textbf{a}, \textbf{b})$ be a $2m$-tuple as defined in the introduction. Let $\mathcal{J}$ denote the set of all subsets $J$ of $[1, m]$ that are of maximal cardinality with respect to the property that $b_j$ for all $j\in \mathcal{J}$ are equal to a fixed integer $b_J$. We note that $\{J: J\in \mathcal{J}\}$ is a partition of $[1, m]$ and that $b_J\not= b_{J'}$ whenever $\{J, J'\}\subseteq \mathcal{J}$. Because $(a_i, b_i)\not= (a_j, b_j)$ whenever $i\not=j$, it follows that if $J\in \mathcal{J}$ then the integers $a_j$ for $j\in J$ are all distinct. Let
 \[
 S_J=\bigcup_{j\in J}\ \{a_j+b_Ji: i \in [0, s-1]\},\ J\in \mathcal{J}.
 \]
Then
 \begin{equation*}
\bigcup_{j=1}^{m}\ \{a_j+b_j(n+i): i\in [0, s-1]\}=\bigcup_{J\in \mathcal{J}}\ (b_Jn+S_J)\ ,\ \forall\ n\in [1, +\infty).\tag{3.1}
  \end{equation*} 
It follows that the collection of sets $AP(\textbf{a}, \textbf{b}; s)$ is a special case of the following more general situation. Let $k\in [1, +\infty)$, let $B=\{b_1,\dots,b_k\}$ be a set of positive integers, and let $\textbf{S}=(S_1,\dots,S_k)$ be a $k$-tuple of finite, nonempty subsets of $[0, +\infty)$. By way of analogy with the expression of the elements of $AP(\textbf{a}, \textbf{b}; s)$ according to (3.1), we will denote by $AP(B, \textbf{S})$ the collection of sets defined by
\[
\Big\{\bigcup_{i=1}^k\ (b_in+S_i): n\in [1, +\infty)\Big\}.
\] 
We are interested in the number of elements of $AP(B, \textbf{S})$ that are sets of quadratic residues or, respectively, quadratic non-residues of a prime $p$, and so if $\varepsilon\in \{-1, 1\}$, we let
\[
 c_\varepsilon(p)=|\{A\in AP(B, \textbf{S})\cap 2^{[1, p-1]}: \chi_p(a)=\varepsilon, \ \textrm{for all}\ a\in A\}|.
 \]
and seek an asymptotic formula for $c_\varepsilon(p)$ as $p\rightarrow +\infty$. N.B. We caution the reader to not confuse this definition of  $c_\varepsilon(p)$ with the definition of this symbol that was given in the statement of Theorem 2.3.  Until further notice, the definition that has just been given will be the one that we will use.   
 
Toward that end, begin by noticing that there is a positive constant $C$, depending only on $B$ and $ \textbf{S}$, such that for all $n\geq C$,
 \begin{equation*}
\textrm{the sets $b_in+S_i, i\in [1, k]$, are pairwise disjoint, and}\tag{3.2}
 \end{equation*} 
\begin{equation*}
\textrm{$\bigcup_{i=1}^k\ (b_in+S_i)$ is uniquely determined by $n$}.\tag{3.3} 
 \end{equation*} 
Because of (3.2) and (3.3), if
\[
\alpha=\sum_i|S_i|\ \textrm{and}\  r(p)=\min_i\left[\frac{p-1-\textrm{max}\ S_i}{b_i}\right],
\]
then the sum
\[
2^{-\alpha}\sum_{x=0}^{r(p)}\ \prod_{i=1}^k\ \prod_{j\in S_i}\ \big(1+\varepsilon\chi_p(b_ix+j)\big)
\]
differs from $c_\varepsilon(p)$ by at most $O(1)$, hence, as per the strategy as outlined in the introduction, this sum can be used to determine the asymptotics of $c_\varepsilon(p)$. 

Apropos of that strategy, let
\[
\mathcal{T}=\bigcup_{i=1}^k\ \{(i, j): j\in S_i\},
\] 
and then rewrite the above sum as
 \begin{equation*}
2^{-\alpha}(1+r(p))+2^{-\alpha}\sum_{\emptyset\not= T\subseteq \mathcal{T}}\ \varepsilon^{|T|}\prod_{i=1}^k\ \chi_p(b_i)^{|\{j: (i, j)\in T\}|}\sum_{x=0}^{r(p)}\ \chi_p\Big(\prod_{(i,j)\in T}\ (x+\bar{b_i}j)\Big),\tag{3.4}
 \end{equation*}
 where $\bar{b_i}$ denotes the inverse of $b_i$ modulo $p$, which clearly exists for all $p$ sufficiently large. Our intent now is to estimate the modulus of the second term in (3.4) by means of Lemma 2.2. This term consists of $2^{\alpha}-1$ summands taken over the nonempty subsets $T$ of $\mathcal{T}$, each summand of which is a product of $2^{-\alpha}$, a coefficient of modulus 1, and an incomplete Weil sum of the form
 \[
 \sum_{x=0}^{r(p)}\ \chi_p\Big(\prod_{(i,j)\in T}\ (x+\bar{b_i}j)\Big).
 \]

Let $\Sigma(p)$ denote the second term of the sum in (3.4). In order to estimate $\Sigma(p)$, we must first remove from it the terms to which Lemma 2.2 cannot be applied. Toward that end, let
\begin{quote}
$E(p)=\{\emptyset\not= T\subseteq \mathcal{T}$: the distinct elements, modulo $p$, in the list $\bar{b_i}j, (i, j)\in T$, each occurs an \emph{even} number of times$\}$.
\end{quote} 
We then split $\Sigma(p)$ into the sum $\Sigma_1(p)$  of terms taken over the elements of $E(p)$ and the sum $\Sigma_2(p)= \Sigma(p)- \Sigma_1(p)$. The sum $\Sigma_2(p)$ has no more than  
$2^{\alpha}-1$ terms each of the form
\[
\pm 2^{-\alpha}\sum_{x=0}^{r(p)} \chi_p\Big(\prod_{(i,j)\in T} (x+\bar{b_i}j)\Big),\ \emptyset \not= T\in 2^{\mathcal{T}}\setminus E(p).
\] 
Since $\emptyset \not= T\notin E(p)$, the polynomial in $x$  in this term at which $\chi_p$ is evaluated can be reduced to a product of at least one and no more than $\alpha$ distinct monic linear factors in $x$ over $Z_p$. Hence by Lemma 2.2,
\[
\Sigma_2(p)=O(\sqrt p\log p)\ \textrm{as}\  p\rightarrow +\infty.
\]
We must now estimate
\[
\Sigma_3(p)=2^{-\alpha}(1+r(p))+\Sigma_1(p),
\]
and, as we shall see, it is precisely this term that will produce the dominant term which determines the asymptotic behavior of $c_\varepsilon(p)$.

Since each element of $E(p)$ has even cardinality,
\[
\Sigma_1(p)=2^{-\alpha}\sum_{T \in E(p)}\ \prod_{i=1}^{k}\chi_p(b_i)^{\{j:(i,j)\in T\}|}\ \sum_{x=0}^{r(p)}\ \chi_p\Big(\prod_{(i,j)\in T}(x+\bar{b_i}j)\Big).
\] 
We now examine the sum over $x\in [0, r(p)]$ on the right-hand side of this equation. Because $T \in E(p)$, each term in this sum is either 0 or 1, and a term is 0 precisely when the value of $x$ in that term agrees with the minimal nonnegative residue mod $p$ of $- \bar{b_i}j$, for some element $(i, j)$ of $T$. However, there are at most $\alpha/2$ of these values at which $x$ can agree for each $T\in E(p)$ and so it follows that $\Sigma_3(p)$ differs by at most $O(1)$ from 
\[
\Sigma_4(p)=2^{-\alpha}(1+r(p))\left(1+\sum_{T\in E(p)}\ \prod_{i=1}^{k}\ \chi_p(b_i)^{|\{j:(i,j)\in T\}|}\right).
\]
Consequently,
\begin{equation*}
\textrm{for all $p$ sufficiently large},\  c_\varepsilon(p)-\Sigma_4(p)=O(\sqrt p\log p),\tag{3.5}
\end{equation*}
and so it suffices to calculate $\Sigma_4(p)$ in order to determine the asymptotics of $c_\varepsilon(p)$.

This calculation requires a careful study of $E(p)$. In order to pin this set down a bit more firmly, we make use of the equivalence relation $\approx$ defined on $\mathcal{T}$ as follows: if $((i, j), (l, m))\in  \mathcal{T}\times \mathcal{T} $ then $(i, j)\approx (l, m)$ if $b_lj=b_im$. For all $p$ sufficiently large, $(i, j)\approx (l, m)$  if and only if $\bar{b_i}j\equiv \bar{b_l}m$ mod $p$, and so if we recall that $\mathcal{E}(A)$ denotes the set of all nonempty subsets of even cardinality of a set $A$, then
\begin{quote}
for all $p$ sufficiently large, $E(p)$ consists of all subsets $T$ of $\mathcal{T}$ such that there exists a nonempty subset $\mathcal{S}$ of equivalence classes of $\approx$ and elements $E_S\in \mathcal{E}(S)$ for $S\in \mathcal{S}$ such that
\begin{equation*}
T=\bigcup_{S\in  \mathcal{S}}\ E_S.\tag{3.6}
\end{equation*}
\end{quote}
In particular, it follows that for all $p$ large enough, $E(p)$ does not depend on $p$, hence from now on, we delete the "$p$" from the notation for this set.

The description of $E$ given by (3.6) mandates that we determine the equivalence classes of the equivalence relation $\approx$. In order to do that in a precise and concise manner, it will be convenient to use the following notation: if $b\in [1, +\infty)$ and $S\subseteq [0, +\infty)$, we let $b^{-1}S$ denote the set of all rational numbers of the form $z/b$, where $z$ is an element of $S$. We next let 
\[
\mathcal{K}=\Big\{\emptyset\not= K\subseteq [1, k]: \bigcap_{i\in K}\ b_i^{-1}S_i\not= \emptyset\Big\}.
\]
If $K\in \mathcal{K}$ then we set
\[
S(K)=\bigcap_{i\in K}\ b_i^{-1}S_i
\]
and, with $\textbf{Q}$ denoting the set of rational numbers,
\[
T(K)=S(K)\cap \Big(\bigcap_{i\in [1, k]\setminus K}\ ( \textbf{Q}\setminus b_i^{-1}S_i)\Big).
\]
 Let
\[
\mathcal{K}_{\max}=\{K\in \mathcal{K}: T(K)\not= \emptyset\}.
\]
 Using the theory of linear Diophantine equations, it is then straightforward to verify that the equivalence classes of $\approx$ consist precisely of all sets of the form
\[
\{(i, tb_i): i\in K\},
\]
where $K\in \mathcal{K}_{\max}$ and $t\in T(K)$. If $\sigma\subseteq K$ and $t\in T(K)$, we also set
\[
E(t, \sigma)=\{(i, tb_i): i\in \sigma\}.
\]

Observe next that if the set
\[
\{\{(i, tb_i): i\in K\}: K\in \mathcal{K}, t\in S(K)\}
\]
is ordered by inclusion then the equivalence classes of $\approx$ are the maximal elements of this set. Hence 
 $T(K)\cap T(K')=\emptyset$ whenever $\{K, K'\}\subseteq \mathcal{K}_{\max}$. Consequently, if $(K, K')\in  \mathcal{K}_{\max}\times \mathcal{K}_{\max},\emptyset \not= \sigma\subseteq K,\emptyset \not= \sigma'\subseteq K', t\in T(K)$, and $t'\in T(K')$, then $E(t, \sigma)$ and $ E(t', \sigma')$ are each contained in distinct equivalence classes of $\approx$ if and only if $t\not= t'$ . It now follows from (3.6) and the structure just obtained for the equivalence classes of $\approx$ that
\begin{quote}
if $T\in E$ then there exists a nonempty subset $\mathcal{S}$ of $\mathcal{K}_{\max}$ , a nonempty subset $\Sigma(S)$ of $\mathcal{E}(S)$ for each $S\in \mathcal{S}$ and a nonempty subset $T(\sigma, S)$ of $T(S)$ for each $\sigma\in \Sigma(S)$ and $S\in \mathcal{S}$ such that 
\begin{equation*}
\textrm{the family of sets}\  \Big\{ T(\sigma, S): \sigma\in \Sigma(S),\  S\in \mathcal{S}\Big\}\  \textrm{is pairwise disjoint, and}\tag{3.7}
\end{equation*}
\begin{equation*}
T=\bigcup_{S\in \mathcal{S}}\ \Big[\bigcup_{\sigma\in \Sigma(S)}\Big(\bigcup_{t\in T(\sigma, S)}\ E(t, \sigma)\Big)\Big].\tag{3.8}
 \end{equation*}
\end{quote}

We have now determined via (3.7) and (3.8) the structure of the elements of $E$ in enough detail for effective use in the calculation of $\Sigma_4(p)$. However, if we already know that $c_\varepsilon(p)=0$, the value of  
 $\Sigma_4(p)$ is obviated in our argument. It would hence be very useful to have a way to mediate between the primes $p$ for which $c_\varepsilon(p)= 0$ and the primes $p$ for which  $c_\varepsilon(p)\not= 0$. We will now define and study a gadget which does that.
\section{The $(B, \textbf{S})$-signature of a prime}
Denote by $\Lambda(\mathcal{K})$ the set
\[
\bigcup_{K\in \mathcal{K}_{\max}} \mathcal{E}(K).
\]
Then $\Lambda(\mathcal{K})$ is empty if and only of every element of  $\mathcal{K}_{\max}$ is a singleton.

Suppose that $\Lambda(\mathcal{K})$ is not empty. We will say that $p$ is an \emph{allowable prime} if no element of $B$ has $p$ as a factor. If  $p$ is an allowable prime, then the $(B, \textbf{S})$-\emph{signature of p} is defined to be the multi-set of $\pm 1$'s given by
\[
\Big\{\chi_p\Big(\prod_{i\in I}\  b_i\Big): I\in \Lambda(\mathcal{K}))\Big\}.
\]

\noindent We declare the signature of $p$ to be \emph{positive} if all of its elements are 1, and \emph{non-positive} otherwise. Let
\begin{quote}
$\Pi_+$ (respectively, $ \Pi_-$) denote the set of all allowable primes $p$ such that the $(B, \textbf{S})$-signature of $p$ is positive (respectively, non-positive).
\end{quote}
We can now prove the following two lemmas: the first records some important information about the signature, and the second implies that we need only calculate $\Sigma_4(p)$ for the primes $p$ in $\Pi_+$.
\begin{lem}
\label{lem2}
$(i)$ The set $\Pi_+$ consists precisely of all allowable primes $p$ for which each of the sets
\begin{equation*}
\{b_i: i\in I\},\  I \in \Lambda(\mathcal{K}) ,\tag{**}
\end{equation*}
is either a set of quadratic residues of $p$ or a set of quadratic non-residues of $p$. In particular, $\Pi_+$ is always an infinite set.

$(ii)$ The set $\Pi_-$ consists precisely of all allowable primes $p$ for which at least one of the sets $(**)$ contains a quadratic residue of $p$ and a quadratic non-residue of $p$, $\Pi_-$ is always either empty or infinite, and $\Pi_-$ is empty if and only if for all $I\in \Lambda(\mathcal{K}),\  \prod_{i\in I}\ b_i$ is a square.
\end{lem} 
\emph{Proof}. Suppose that $p$ is an allowable prime such that each of the sets (**) is either a set of quadratic residues of $p$ or a set of quadratic non-residues of $p$. Then
\[
\chi_p\Big(\prod_{i\in I}\ b_i\Big)=1
\]
whenever $I\in  \Lambda(\mathcal{K})$ because $|I|$ is even, i.e., $p\in \Pi_+$. On the other hand, let $p\in \Pi_+$ and let $I=\{i_1,\dots,i_n\} \in  \Lambda(\mathcal{K})$. Then because $p\in \Pi_+$,
\[
\chi_p(b_{i_j}b_{i_{j+1}})=1,\ j\in [1, n-1],
\]
and these equations imply that $\{b_i:i\in I\}$ is either a set of quadratic residues of $p$ or a set of quadratic non-residues of $p$. This verifies the first statement in ($i$), and the second statement follows from the fact that there are infinitely many primes $p$ such that $B$ is a set of quadratic residues of $p$.

Statement ($ii$) of the lemma follows from ($i$), the definition of $\Pi_-$, and the fact that a positive integer is a quadratic residue of all but finitely many primes if and only if it is a square. $\ \Box$
\begin{lem}
\label{lem3}
 If $p\in \Pi_-$ then $c_\varepsilon(p)=0$.
 \end{lem}
\emph{Proof}. If $p\in \Pi_-$ then there is an $I\in  \Lambda(\mathcal{K})$ such that
\[
\chi_p\Big(\prod_{i\in I} b_i\Big)=-1.
\]
Because $I$ is nonempty and of even cardinality, there exists $\{m, n\}\subseteq I$ such that
\begin{equation*}
\chi_p(b_mb_n)=-1.\tag{4.1}
\end{equation*}
Because $\{m, n\}$ is contained in an element of  $\mathcal{K}_{\max}$, it follows that $b_m^{-1}S_m\cap b_n^{-1}S_n\not= \emptyset$, and so we find a non-negative rational number $r$ such that
\begin{equation*}
rb_m\in S_m\ \textrm{and}\ rb_n\in S_n.\tag{4.2}
\end{equation*}
By way of contradiction, suppose that $c_\varepsilon(p)\not= 0$. Then there exists a $z\in [1, +\infty)$ such that $b_mz+S_m$ and $b_nz+S_n$ are both contained in $[1, p-1]$ and 
\begin{equation*}
\chi_p(b_mz+u)=\chi_p(b_nz+v),\ \textrm{for all}\  u\in S_m\ \textrm{and for all}\ v\in S_n.\tag{4.3}
\end{equation*}

If $d$ is the greatest common divisor of $b_m$ and $b_n$ then there is a non-negative integer $t$ such that $r=t/d$. Hence by (4.2) and (4.3),
\begin{eqnarray*}
\chi_p(b_m/d)\chi_p(dz+t)&=&\chi_p(b_mz+rb_m)\\
&=&\chi_p(b_nz+rb_n)\\
&=&\chi_p(b_n/d)\chi_p(dz+t).
\end{eqnarray*}
However, $dz+t\in [1, p-1]$ and so $\chi_p(dz+t)\not= 0$, hence
\[
\chi_p(b_m/d)=\chi_p(b_n/d).
\]
But then
\[
\chi_p(b_mb_n)=\chi_p(d^2)\chi_p(b_m/d)\chi_p(b_n/d)=1,
\]
contrary to (4.1).   $\ \Box$
\section{The calculation of $\Sigma_4(p)$}
In this section we calculate the sum $\Sigma_4(p)$ that arose from the work of section 3.
By virtue of Lemma 4.2, we need only calculate $\Sigma_4(p)$ for $p\in \Pi_+$, hence let $p$ be an allowable prime for which
\begin{equation*}
\chi_p\Big(\prod_{i\in I} b_i\Big)=1,\ \textrm{for all}\ I\in  \Lambda(\mathcal{K}).\tag{5.1}
\end{equation*}

We first recall that
\begin{equation*}
\Sigma_4(p)=2^{-\alpha}(1+r(p))\left(1+\sum_{T\in E}\ \prod_{i=1}^{k}\ \chi_p(b_i)^{|\{j: (i, j)\in T\}|}\right),\tag{5.2}
\end{equation*}
and so we must evaluate the products over $T\in E$ which determine the summands of the third factor on the right-hand side of (5.2). Toward that end, let $T\in E$ and find a nonempty subset $ \mathcal{S}$ of $\mathcal{K}_{\max}$, a nonempty subset $\Sigma(S)$ of $\mathcal{E}(S)$ for each $S\in \mathcal{S}$ and a nonempty subset $T(\sigma, S)$ of $T(S)$ for each $\sigma\in \Sigma(S)$ and $S\in \mathcal{S}$ such that  (3.7) holds and $T$ satisfies (3.8).Then
\begin{eqnarray}
\{j: (i, j)\in T\}&=&\bigcup_{S\in \mathcal{S}}\ \Big[\bigcup_{\sigma\in \Sigma(S)}\Big(\bigcup_{t\in T(\sigma, S)} \{j: (i, j)\in E(t, \sigma)\}\Big)\Big]\label{2.13}\\
&=&\bigcup_{S\in \mathcal{S}}\ \Big(\bigcup_{\sigma\in \Sigma(S): i\in \sigma}\{tb_i: t\in T(\sigma, S)\}\Big).
\end{eqnarray}
 It follows from (3.7) that the union (2) is pairwise disjoint. Hence
\[
|\{j: (i, j)\in T\}|=\sum_{S\in  \mathcal{S}}\sum_{\sigma\in \Sigma(S): i\in \sigma}\ |T(\sigma, S)|.
\]
Thus from this equation and (5.1) we find that
\begin{eqnarray*}
\prod_{i=1}^{k}\ \chi_p(b_i)^{|\{j: (i, j)\in T\}|}&=&\prod_{i\in \cup_{S\in \mathcal{S}}\cup_{\sigma\in \Sigma(S)}\ \sigma}\ \chi_p(b_i)^{\sum_{S\in  \mathcal{S}}\sum_{\sigma\in \Sigma(S): i\in \sigma} |T(\sigma, S)|}\\
&=&\prod_{S\in  \mathcal{S}}\Big(\prod_{\sigma\in \Sigma(S)}\ \Big(\chi_p\Big(\prod_{i\in \sigma} b_i\Big)\Big)^{|T(\sigma, S)|}\Big)\\
&=&1.
\end{eqnarray*}
Hence
\begin{equation*}
\sum_{T\in E}\ \prod_{i=1}^{k}\ \chi_p(b_i)^{|\{j: (i, j)\in T\}|}=|E|,\tag{5.3}
\]
and so we must count the elements of $E$. In order to do that, note first that the pairwise disjoint decomposition (3.6) of an element $T$ of $E$ is uniquely determined by $T$, and, obviously, uniquely determines $T$. Hence if $\mathcal{D}$ denotes the set of all equivalence classes of $\approx$ of cardinality at least 2 then
\begin{eqnarray*}
|E|&=&\sum_{\emptyset\not= \mathcal{S}\subseteq \mathcal{D}}\ \prod_{S\in \mathcal{S}}\ |\mathcal{E}(S)|\\
&=&-1+\prod_{D\in \mathcal{D}} (1+|\mathcal{E}(D)|)\\
&=&-1+ \prod_{D\in \mathcal{D}} 2^{|D|-1}\\
&=&-1+2^{-|\mathcal{D}|}\cdot 2^{\sum_{D\in \mathcal{D}}|D|}. 
\end{eqnarray*}
However, $\mathcal{D}$ consists of all sets of the form
\[
\{(i, tb_i): i\in K\}
\]
where $K\in \mathcal{K}_{\max}, |K|\geq 2$, and $t\in T(K)$. Hence  
\[
|\mathcal{D}|=\sum_{K\in  \mathcal{K}_{\max}: |K|\geq 2}\ |T(K)|,
\]
\[
\sum_{D\in \mathcal{D}}|D|=\sum_{K\in  \mathcal{K}_{\max}: |K|\geq 2}\ |K||T(K)|,
\]
and so if we set
\[
e=\sum_{K\in  \mathcal{K}_{\max}}\ |T(K)|(|K|-1),
\] 
then
\begin{equation*}
|E|=2^e-1.\tag{5.4}
\]

\noindent Equations (5.2), (5.3), and (5.4) now imply
\begin{lem}
\label{lem4}
If
\[
\alpha=\sum_i|S_i|,\ e=\sum_{K\in  \mathcal{K}_{\max}}\ |T(K)|(|K|-1),\ \textrm{and}\  r(p)=\min_i\left[\frac{p-1-\max S_i}{b_i}\right],
\] 
then
\[
\Sigma_4(p)=2^{e-\alpha}(1+r(p)),\ \textrm{for all}\  p\in \Pi_+.
\]
\end{lem}
\section{The asymptotic behavior of $c_\varepsilon(p)$}
With Lemmas 4.1, 4.2, and 5.1 now in hand, we can prove the following theorem, in which the asymptotic behavior of $c_\varepsilon(p)$ is determined.
\begin{thm}
\label{thm5}
Let $\varepsilon\in \{-1, 1\}, k\in [1, +\infty)$, and let $B=\{b_1,\dots,b_k\}$ be a set of positive integers and $\textbf{S}=(S_1,\dots,S_k)$ a $k$-tuple of finite, nonempty subsets of $[0, +\infty)$. If $ \mathcal{K}_{\max}$ is the set of  subsets of $[1, k]$ as determined in section 3 by $B$ and  $\textbf{S}$, let 
\[
\Lambda(\mathcal{K})=\bigcup_{K\in \mathcal{K}_{\max}} \mathcal{E}(K),
\]
\[
\alpha=\sum_i|S_i|,\ b=\max_i\{b_i\},\  e=\sum_{K\in  \mathcal{K}_{\max}}\ |T(K)|(|K|-1),\ and
\]
\[
 c_\varepsilon(p)=|\{A\in AP(B, \textbf{S})\cap 2^{[1, p-1]}: \chi_p(a)=\varepsilon, \ \textrm{for all}\ a\in A\}|.
\]
$(i)$ If $\Lambda(\mathcal{K})$ is empty then
\[
c_\varepsilon(p)\sim (b\cdot 2^\alpha)^{-1}p\ as\ p\rightarrow +\infty.
\]

\noindent$(ii)$ If $\Lambda(\mathcal{K})$ is not empty then

$(a)$ the parameter $e$ is positive;

$(b)$ if $\prod_{i\in I} b_i$ is a square for all $I\in \Lambda(\mathcal{K})$ then
\[
c_\varepsilon(p)\sim (b\cdot 2^{\alpha-e})^{-1}p\ as\ p\rightarrow +\infty;
\]

$(c)$ if there exists $I\in \Lambda(\mathcal{K})$ such that $\prod_{i\in I} b_i$ is not a square then

$(\alpha)$ the set $\Pi_+$ of primes with positive $(B, \textbf{S})$-signature and the set $\Pi_-$ of primes with non-positive $(B, \textbf{S})$-signature are both infinite,

$(\beta)$ $c_\varepsilon(p)=0$ for all $p$ in $\Pi_-$, and

$(\gamma)$ as $p\rightarrow +\infty$ inside $\Pi_+$,
\[
c_\varepsilon(p)\sim (b\cdot2^{\alpha-e})^{-1}p\ .
\]
 \end{thm}
 \emph{Proof}. If $\Lambda(\mathcal{K})$ is empty then every element of $\mathcal{K}_{\max}$ is a singleton set, hence all of the equivalence classes of the equivalence relation $\approx$ defined above on $\mathcal{T}$ by the set $B$ are singletons. It follows that the set $E$ which is summed over in (5.2) is empty and so  
\begin{equation*}
\Sigma_4(p)=2^{-\alpha}(1+r(p)),\ \textrm{for all}\ p\ \textrm{sufficiently large}.\tag{6.1}
\]
Upon recalling that
\[
r(p)=\min_i\left[\frac{p-1-\textrm{max}\ S_i}{b_i}\right],
\]
the conclusion of $(i)$ is an immediate consequence of (3.5) and (6.1).

Suppose that $\Lambda(\mathcal{K})$ is not empty. Conclusion $(a)$ is obvious. If $\prod_{i\in I} b_i$ is a square for all $I\in \Lambda(\mathcal{K})$ then it follows from its definition that $\Pi_+$ contains all but finitely many primes, and so $(b)$ is an immediate consequence of (3.5) and Lemma 5.1. On the other hand, if there exists $I\in \Lambda(\mathcal{K})$ such that $\prod_{i\in I} b_i$ is not a square then $(\alpha)$ follows from Lemma 4.1, $(\beta)$ follows from Lemma 4.2, and $(\gamma)$  is an immediate consequence of (3.5) and Lemma 5.1.  $\   \Box$  

Theorem 6.1 shows that the elements of $\Lambda(\mathcal{K})$ contribute to the formation of quadratic residues and non-residues inside $AP(B, \textbf{S})$. If no such elements exist then $c_\varepsilon(p)$ has the expected minimal asymptotic approximation $(b\cdot 2^{\alpha})^{-1}p$ as $p\rightarrow +\infty$. In the presence of elements of  $\Lambda(\mathcal{K})$, the parameter $e$ is positive, the  asymptotic size of $c_\varepsilon(p)$ is increased by a factor of $2^e$, and whenever $\Pi_-$ is empty, $c_\varepsilon(p)$ is asymptotic to $(b\cdot 2^{\alpha-e})^{-1}p$ as $p\rightarrow +\infty$. However, the most interesting behavior occurs when $\Pi_-$ is not empty; in that case, as $p\rightarrow +\infty, c_\varepsilon(p)$ asymptotically oscillates infinitely often between 0 and $(b\cdot 2^{\alpha-e})^{-1}p$. 

\emph{Remarks}. (1) We note that the proof of Theorem 6.1 yields the following error estimates for the asymptotic approximations in that theorem: if the hypothesis of $(i)$ is satisfied, then for all $p$ sufficiently large,
\[
|c_\varepsilon(p)-(b\cdot 2^{-\alpha})^{-1}p|\leq (1+2\alpha)\sqrt p\log p,
\]
and if the hypothesis of $(ii)$ and $(b)$ (respectively, $(c)$) is satisfied, then for all $p$ sufficiently large, (respectively, for all $p$ sufficiently large inside $\Pi_+$),
\[
|c_\varepsilon(p)-(b\cdot 2^{\alpha-e})^{-1}p|\leq (1+2\alpha)\sqrt p\log p. 
\]

(2)  In [7], Davenport in fact considered and solved a more general problem than the one that is discussed in section 1. As we pointed out at the beginning of section 2, he showed in [7, Corollary of Theorem 5] that if $\eta$ is a choice of signs for $[0, s-1]$ then the cardinality of the set
\[
\{\{x+i: i\in [0, s-1]\}\in AP(0, 1;s)\cap 2^{[1, p-1]}: \chi_p(x+i)=\eta(i), i\in [0, s-1]\}
\]
is asymptotic to $2^{-s}p$ as $p\rightarrow +\infty$, i.e., he asymptotically enumerated the elements of  $AP(0, 1;s)\cap 2^{[1, p-1]}$ which exhibit a fixed but arbitrary pattern of quadratic residues and non-residues of $p$. When the choice of signs is either identically 1 or identically $-1$, we recover the results that are discussed in section 1.

An analog of this more general problem can also be formulated in the context of our work here. For $k\in [1, +\infty)$ and $\textbf{S}=(S_1,\dots, S_k)$ a $k$-tuple of nonempty subsets of $[0, +\infty)$, let $\eta_i$ be a choice of signs for $S_i, i\in [1, k]$. Let $B=\{b_1,\dots, b_k\}$ be a subset of  $[1, +\infty)$. Setting $\eta=(\eta_1,\dots, \eta_k)$, we let $c_{\eta}(p)$ denote the cardinality of the set
 \[
\Big\{\bigcup_{i=1}^k\  (b_ix+S_i)\in AP(B, \textbf{S})\cap 2^{[1, p-1]}: \chi_p(b_ix+j)=\eta_i(j), \forall\ i \in [1, k], \forall\  j\in S_i\Big\},
 \]
and we then consider the problem of determining the asymptotic behavior of $c_{\eta}(p)$ as $p\rightarrow +\infty$. One easily deduces from the arguments of section 3 that if $\alpha=\sum_i|S_i|$, $r(p)=\min_i\{(p-1-\max S_i)/b_i\}$, $E$ is the set defined by (3.6), and $\Sigma_{\eta}(p)$ is the sum
\begin{equation*}
2^{-\alpha}(1+r(p))\left(1+\sum_{T\in E}\ \prod_{(i, j)\in T}\ \eta_i(j)\ \prod_{i=1}^{k}\ \chi_p(b_i)^{|\{j:(i,j)\in T\}|}\right),\tag{6.2}
\end{equation*}
then
\begin{equation*}
c_{\eta}(p)-\Sigma_{\eta}(p)=O(\sqrt p\log p),\ \textrm{as}\  p\rightarrow +\infty,\tag{6.3}
\end{equation*}
and so we conclude that whenever the set $\Lambda(\mathcal{K})$ defined in Theorem 6.1 is empty and $b=\max_i\{b_i\}$ then
\[ 
c_{\eta}(p)\sim (b\cdot 2^{\alpha})^{-1}p\ \textrm{as}\  p\rightarrow +\infty,
\]
i.e., Theorem 6.1($i$) remains valid. On the other hand, if $\Lambda(\mathcal{K})$ is not empty then the constant $e$ defined in Theorem 6.1 is positive, and from our calculation of $|E|$ in the proof of Lemma 5.1, it follows that
\[
\left|\Sigma_{\eta}(p)\right|\leq 2^{e-\alpha}(1+r(p)),
\]
but because of the sign $\prod_{(i, j)\in T} \eta_i(j)$ in the terms of $\Sigma_{\eta}(p)$, we have been unable to finish the calculation of $\Sigma_{\eta}(p)$ in this case.  Hence the ideas of this paper  apparently provide no further insight into the asymptotic behavior of $c_{\eta}(p)$.

However, there \emph{is} a special case of this problem for which our methods can be pushed through to give a nontrivial generalization of Theorem 6.1. Let $\eta$ now denote a choice of signs for $[1, k]$ and take the choice of signs $\eta_i$ on $S_i$ to be identically $\eta(i)$, i.e., we take the choice of signs on $S_i$ to be \emph{constant} for each $i\in [1, k]$. In a conflation of notation that we hope will not be confusing, we also let $\eta$ denote the corresponding $k$-tuple of the choice of signs that we have just defined for the $S_i$'s and then observe that the associated sum (6.2) becomes 
\[ 
 2^{-\alpha}(1+r(p))\left(1+\sum_{T\in E}\  \prod_{i=1}^{k}\ (\eta(i)\chi_p(b_i))^{|\{j:(i,j)\in T\}|}\right).
 \]
 When $\Lambda(\mathcal{K})$ is not empty, we now define the $(B, \textbf{S}, \eta)-signature$ of an allowable prime $p$ to be the multi-set of $\pm1$'s given by
 \[
 \Big\{\prod_{i\in I}\eta(i)\chi_p(b_i): I\in \Lambda(\mathcal{K})\Big\},
 \]
and then define the sets of primes $\Pi_+$ and $\Pi_-$ as it was done previously in Section 4, with the $(B, \textbf{S}, \eta)$-signature in place of the $(B, \textbf{S})$-signature. The set $\Pi_-$ is hence the complement of $\Pi_+$ in the set of allowable primes, and the proof of Lemmas 4.1 and 4.2 is easily modified to show that $\Pi_+$  consists of all allowable primes $p$ such that the function $i\rightarrow \eta(i)\chi_p(b_i)$ defined on $[1, k]$ is constant on each set in $\Lambda(\mathcal{K})$, and also that
\begin{equation*}
c_{\eta}(p)=0,\ \textrm{for all}\ p\in \Pi_-.\tag{6.4}
\]
If $ \Lambda(\mathcal{K})$ is nonempty then we also calculate as before that
\begin{equation*}
\Sigma_{\eta}(p)=2^{e-\alpha}(1+r(p)),\ \textrm{for all}\  p\in \Pi_+.\tag{6.5}
\end{equation*} 
As we already noted, Theorem 6.1($i$) remains valid, and by virtue of (6.3)-(6.5), the validity of the rest of Theorem 6.1 in this more general context depends only on the structure of  $\Pi_+$ and  $\Pi_-$. If $\Pi_+$ contains all but finitely many primes then the conclusions of Theorem 6.1($ii$) ($a$) and ($b$) hold, if $\Pi_+$ and $\Pi_-$ are both infinite then the conclusions of Theorem 6.1($ii$) ($\alpha$), ($\beta$), and ($\gamma$) hold, and if  $\Pi_+$ is finite then $c_{\eta}(p)=0$ for all $p$ sufficiently large.

Thus our attention is focused on the following interesting problem: given $k\in [1, +\infty)$, a choice of signs $\eta$ on $[1, k]$, a subset $B=\{b_1,\dots,b_k\}$ of $[1, +\infty)$, and a $k$-tuple $\textbf{S}$ of nonempty subsets of $[0, +\infty)$, characterize when the corresponding set $\Pi_+$ either contains all but finitely many primes, is infinite, or is finite. A solution of this problem will involve rather delicate combinatorial relationships between $\eta$, the elements of $\Lambda(\mathcal{K})$, and the prime factorization of the elements of $B$; one such solution can be found by using the results and methods of [16] and [17]. In particular, one can prove that $\Pi_+$ is always either empty or infinite and $\Pi_+$ contains all but finitely many primes if and only if for all $I\in \Lambda(\mathcal{K})$, $\eta$ is constant on $I$ and $\prod_{i\in I} b_i$ is a square . Because it would take us too far afield at this point, we leave the verification of these facts and the other details to the interested reader. 

\section{Discussion of the asymptotics associated with $AP(\textbf{a}, \textbf{b}; s)$}
Theorem 6.1 will now be applied to the situation of primary interest to us here, namely to the family of sets  $AP(\textbf{a}, \textbf{b}; s)$ as defined in section 1 by a $2m$-tuple $(\textbf{a}, \textbf{b})$. In this case, the decomposition (3.1) of the sets in $AP(\textbf{a}, \textbf{b}; s)$ shows that there is a set $B=\{b_1,\dots,b_k\}$ of positive integers, a $k$-tuple $(m_1,\dots,m_k)$ of positive integers such that $m=\sum_i m_i$, and sets
\[
A_i=\{a_{i1},\dots,a_{im_i}\}
\]
of non-negative integers such that if we let
\begin{equation*}
S_i=\bigcup_{j=1}^{m_i}\ \{a_{ij}+b_il: l\in [0, s-1]\},\  i\in [1, k],\tag{7.1}
\]  
and set
\[
\textbf{S}=(S_1,\dots,S_k)
\]
then
\[
AP(\textbf{a}, \textbf{b}; s)= AP(B, \textbf{S}).
\]
After letting $Q_i$ denote the set of rational numbers obtained when the elements of the set $A_i$ are divided by $b_i$, it follows that
\[
b_i^{-1}S_i=\bigcup_{q\in Q_i}\ \{q+j: j\in [0, s-1]\},\ i\in [1, k].
\]
These sets then determine the subsets of $[1, k]$ that constitute 
\[
\mathcal{K}=\{\emptyset\not= K\subseteq [1, k]: \bigcap_{i\in K}\ b_i^{-1}S_i\not= \emptyset\}\}
\]
and hence also the elements of $\mathcal{K}_{\max}$, according to the recipe given in section 3. 
The sets in $\mathcal{K}_{\max}$, together with the parameters
\[
\alpha=\sum_i|S_i|,\ b=\max_i\{b_i\},\  \textrm{and}\  e=\sum_{K\in  \mathcal{K}_{\max}}\ |T(K)|(|K|-1),
\]
when used as specified in Theorem 6.1, then determine precisely the asymptotic behavior of the sequence $c_\varepsilon(p)$ that is defined upon replacement of $AP(B, \textbf{S})$ by $AP(\textbf{a}, \textbf{b}; s)$ in the statement of Theorem 6.1. In particular, $\Lambda(\mathcal{K})$ is empty if and only if 
 \begin{quote}
 $(7.2)$  if $(i, j)\in [1, k]\times [1, k]$ with $i\not= j$ and $(a, a')\in A_i\times A_j$, then either $b_ib_j$ does not divide $a'b_i-ab_j$ or  $b_ib_j$ divides $a'b_i-ab_j$ with a quotient that exceeds $s-1$ in modulus.
 \end{quote}
Hence the conclusion of statement $(i)$ of Theorem 6.1 holds for $AP(\textbf{a}, \textbf{b}; s)$ when condition (7.2) is satisfied, while the conclusions of statement $(ii)$ of Theorem 6.1 hold for $AP(\textbf{a}, \textbf{b}; s)$ whenever condition (7.2) is not satisfied. In section 8, we will show, among other things, that for $m\in [2, +\infty)$ and for each of the hypotheses in the statement of Theorem 6.1, there exists infinitely many $2m$-tuples $(\textbf{a}, \textbf{b})$ which satisfy that hypothesis. 

In addition to the parameter $e$, the sum $\alpha$ of the cardinalities of the sets $S_i$ which are defined by (7.1) is an important parameter in the coefficient of $p$ in the asymptotic formula of $c_\varepsilon(p)$. One can use the principle of inclusion and exclusion to calculate $\alpha$, but there is an alternative calculation of the cardinality of these sets which in practice is often more tractable than the calculation which uses inclusion and exclusion. Because of its relevance to our discussion here, we will now carry it out.

The calculation is based on the concept of  what we will call an \emph{overlap diagram} (overlap diagrams will also be used in some calculations that we will perform in section 8). In order to define this diagram, let $(n, s)\in [1,+\infty)\times [1,+\infty)$ and let $ \textbf{g}=(g(1),\dots,g(n))$ be an $n$-tuple of positive integers. We use \textbf{g} to construct the following array of points. In the plane, place $s$ points horizontally one unit apart, and label the $j$-th point as $(1, j-1)$ for each $j\in [1, s]$. This is \emph{row $1$}. Suppose that row $i$ has been defined. One unit vertically down and $g(i)$ units horizontally to the right of the first point in row $i$, place $s$ points horizontally one unit apart, and label the $j$-th point as $(i+1, j-1)$ for each $j\in [1, s]$. This is $\emph{row}\ i+1$.
The array of points so formed by these $n+1$ rows is called the \emph{overlap diagram of} \textbf{g}, the sequence \textbf{g} is called the \emph{gap sequence} of the overlap diagram, and a nonempty set that is formed by the intersection of the diagram with a vertical line is called a \emph{column} of the diagram. N.B. We do not distinguish between the different possible positions in the plane which the overlap diagram may occupy. A typical example with $n=3, s=8$, and gap sequence (3, 2, 2) looks like
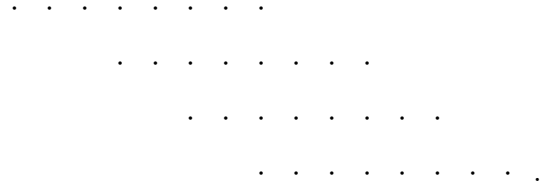
\vspace{0.5cm}
\begin{center}
\begin{tabular}{ccccccccccccccc}
 $\cdot$&$\cdot$&$\cdot$&$\cdot$&$\cdot$&$\cdot$&$\cdot$&$\cdot$&&&&&&&\\

&&&$\cdot$&$\cdot$&$\cdot$&$\cdot$&$\cdot$&$\cdot$&$\cdot$&$\cdot$&&\\
&&&&&$\cdot$&$\cdot$&$\cdot$&$\cdot$&$\cdot$&$\cdot$&$\cdot$&$\cdot$\\
&&&&&&&$\cdot$&$\cdot$&$\cdot$&$\cdot$&$\cdot$&$\cdot$&$\cdot$&$\cdot$\     \ .\\
\end{tabular}
\end{center}
\vspace{0.5cm} 

 We need to describe how and where rows overlap in an overlap diagram. Begin by first noticing that if $(g(1),\dots,g(n))$ is the gap sequence, then row $i$ overlaps row $j$ for $i<j$ if and only if
 \[
 \sum_{r=i}^{j-1} g(r)\leq s-1;
 \]                                                         
in particular, row $i$ overlaps row $i+1$ if and only if $g(i)\leq s-1$. Now let $\mathcal{G}$ denote the set of all subsets $G$ of $[1, n]$ such that $G$ is a nonempty set of consecutive integers maximal with respect to the property that $g(i)\leq s-1$ for all $i\in G$. If $\mathcal{G}$ is empty then $g(i)\geq s$ for all $i\in [1, n]$, and so there is no overlap of rows in the diagram. Otherwise there exists $m\in[1, 1+[n/2]]$ and strictly increasing sequences $(l_1,\dots,l_m)$ and $(M_1,\dots,M_m)$ of positive integers, uniquely determined by the gap sequence of the diagram, such that $l_i\leq M_i$ for all $i\in [1, m], 1+M_i\leq l_{i+1}$ if $i\in [1, m-1]$, and
\[
\mathcal{G}=\{[l_i, M_i]: i\in [1, m]\}.
\]
In fact, $l_{i+1}> 1+M_i$ if $i\in [1, m-1]$, lest the maximality of the elements of $\mathcal{G}$ be violated. It follows that the intervals of integers $[l_i, 1+M_i], i\in [1, m]$, are pairwise disjoint. 

The set $\mathcal{G}$ can now be used to locate the overlap between rows in the overlap diagram like so: for $i\in [1, m]$, let
\[
B_i=[l_i, 1+M_i],
\]
and set
\[
\mathcal{B}_i=\textrm{the set of all points in the overlap diagram whose labels are in}\ B_i\times [0, s-1].
\]
We refer to $\mathcal{B}_i$ as the \emph{i-th block} of the overlap diagram, to the interval of integers $B_i$ as the \emph{support of} $\mathcal{B}_i$, and to the sequence $(g(j): j\in [l_i, M_i])$ as the \emph{gap sequence} of $\mathcal{B}_i$. Thus the blocks of the diagram are precisely the regions in the diagram in which rows overlap. 

Now let $(b, t)\in [1, +\infty)\times [1, +\infty)$ and let $\{a_1,\dots,a_t\}$ be a set of non-negative integers. If we set $r_i=a_i/b$ then with no loss of generality , we assume that the $a_i$'s are indexed so that $r_i<r_{i+1}$ for each $i\in [1, t-1]$. We will use overlap diagrams to calculate
\[
\mu=\Big|\bigcup_{j=1}^t\ \{a_j+bi: i\in [0, s-1]\}\Big|.
\]

To that end, then, we define the equivalence relation $\approx_b$ on $[1, t]\times [0, s-1]$ by
\[
(i, j)\approx_b (l, m) \ \textrm{if}\ a_i-a_l=b(m-j);
\] 
$\mu$ is the total number of equivalence classes of $\approx_b$. Next, let
\[
Q=\{(i, j)\in [1, t]\times [1, t]: i\not= j\ \textrm{and}\ b \ \textrm{divides}\ a_i-a_j\ \textrm{with quotient}\ q(i, j)\}.
\]
The quotients $q(i, j)$ are all nonzero because $a_i\not= a_j$ for $i\not=j$.  If $\pi$ denotes the canonical projection of $[1, t]\times [1, t]$ onto its left factor, then on $\pi(Q)$ we consider the equivalence relation defined by $i\sim j$ if $i=j$ or $(i, j)\in Q$, and we let $\{E_1,\dots,E_v\}$ denote the equivalence classes of $\sim$.  We have that
\[
r_i-r_j=q(i, j),\ \textrm{for all}\ (i, j)\in Q,
\]
and so each set of rational numbers $\{r_i: i\in E_n\}, n\in [1, v]$, is linearly ordered accordingly. The elements of this set are listed in increasing order,  and so we let $q_n(i)$ denote the  quotient that equals the positive difference between element $i$ and element $i+1$ on that list, $i\in  [1, |E_n|-1]$,  and then let $\mathcal{D}_n$ denote the overlap diagram of the $(|E_n|-1)$-tuple  $(q_n(i): i\in [1, |E_n|-1])$. 

Observe now that the equations
\[
r_i-r_l=m-j
\]
define an equivalence relation on $E_n\times [0, s-1]$ such that the number of columns $c(n)$ of $\mathcal{D}_n$  counts the equivalence classes of this equivalence relation. However, we have that
\[
r_i-r_l=m-j\ \textrm{if and only if}\ a_i-a_l=b(m-j),
\]
and so it follows that
\begin{quote}
(7.2)\ \    $c(n)$ is the number of equivalence classes of $\approx_b$ that are determined by the elements of $E_n\times [0, s-1]$.
\end{quote}  
Also, if $i\not= j$ then the set of equivalence classes of $\approx_b$ that are determined by the elements of $E_i\times [0, s-1]$ is disjoint from the set of equivalence classes of  $\approx_b$ that are determined by the elements of  $E_j\times [0, s-1]$. We conclude from (7.2) that
\begin{quote}
(7.3)\ \  $\sum_n c(n)$ is the number of equivalence classes of $\approx_b$ that are determined by the elements of $\pi(Q)\times [0, s-1]$.
\end{quote}  
On the other hand, it follows from the definition of $Q$ that the equivalence classes of $\approx_b$ that are determined by elements of $([1, t]\setminus \pi(Q))\times [0, s-1]$ are singletons  and this set of equivalence classes is disjoint from the set of equivalence classes that are  determined by the elements of $\pi(Q)\times [0, s-1]$. 
Therefore, because of (7.3),
\begin{equation*}
\mu=s(t-|\pi(Q)|)+\sum_n c(n).\tag{7.4}
\end{equation*}

The next step is to use the block structure of $\mathcal{D}_n$ to calculate $\sum_n c(n)$. Let
\[
\mathcal{O}=\{n\in [1, v]:\ \textrm{there exists}\ (i, j)\in E_n\times E_n\ \textrm{such that}\ |q(i, j)|\leq s-1\}.
\]
If $n \in [1, v]\setminus \mathcal{O}$ then there is no overlap of rows in $\mathcal{D}_n$, and so
\begin{equation*}
\textrm{if}\ n\in [1, v]\setminus \mathcal{O}\ \textrm{then}\ c(n)=s|E_n|.\tag{7.5}
\end{equation*}
On the other hand, if $n\in \mathcal{O}$ then there are rows of $\mathcal{D}_n$ which overlap to form blocks $\mathcal{B}_r(n)$ for $r\in [1, m(n)]$, say, with corresponding supports $B_r(n), r\in [1, m(n)]$. 
The gap sequence of $\mathcal{B}_r(n)$ can be indexed as 
\[
(q_{nr}(i): i\in [1, |B_r(n)|-1])
\]
for a suitable subsequence of the gap sequence of $\mathcal{D}_n$. A moment's reflection now reveals that
\begin{equation*}
\textrm{the number of columns of}\  \mathcal{B}_r(n)\  \textrm{is}\  s+\sum_{i=1}^{|B_r(n)|-1} q_{nr}(i).\tag{7.6}
\end{equation*} 
After observing that the columns of $\mathcal{D}_n$ that are not contained in a block  are in one-to-one correspondence with the elements of $[1, |E_n|]\times [0, s-1]$ outside of
\[
\bigcup_{r=1}^{m(n)}\ (B_r(n)\times [0, s-1]),
\]
and also recalling that the supports $B_r(n)$ are pairwise disjoint, we conclude from (7.6) that
\begin{equation*}
\textrm{if}\ n\in \mathcal{O}\ \textrm{then}\ c(n)=s(|E_n|+m(n)-\sum_{r=1}^{m(n)}\ |B_r(n)|)+\sum_{r=1}^{m(n)}\ \sum_{i=1}^{|B_r(n)|-1}q_{nr}(i).\tag{7.7}
\end{equation*}
From (7.5) and (7.7), we conclude that
\begin{equation*}
\sum_n c(n)=s|\pi(Q)|+\sum_{n\in \mathcal{O}}\ \sum_{r=1}^{m(n)}\ \sum_{i=1}^{|B_r(n)|-1} (q_{nr}(i)-s).\tag{7.8}
\end{equation*}

It now follows from (7.4) and  (7.8) that
\begin{equation*}
\mu=st-\sum_{n\in \mathcal{O}}\ \sum_{r=1}^{m(n)}\ \sum_{i=1}^{|B_r(n)|-1} (s-q_{nr}(i)).\tag{7.9}
\end{equation*}

We will now reformulate (7.9) so that a calculation of $\mu$ can be carried out using information that comes in a straightforward manner directly from the set $\{r_1,\dots,r_t\}$,  thereby dispensing with the use of overlap diagrams.

Begin by letting $\mathcal{R}$ denote the set of all subsets $R$ of $\{r_1,\dots,r_t\}$ such that $|R|\geq 2$ and $R$ is maximal with respect to the property that $u-v$ is an integer for all $(u, v)\in R\times R$. If $R\in \mathcal{R}$ then we linearly order the elements of $R$ and let $D(R)$ denote the $(|R|-1)$-tuple of positive integers whose coordinates are the distances between consecutive elements of $R$ in this ordering. Now observe that there is a natural bijective correspondence between the equivalence classes $E_n$ and the elements of $\mathcal{R}$ such that when $E_n$ and $R$ correspond, the gap sequence of the diagram $\mathcal{D}_n$ is given by the $(|R|-1)$-tuple $D(R)$. Consequently, (7.9) immediately implies the following proposition:  
\begin{prp}
\label{prp1}
If for each $R\in \mathcal{R}$  we let $M_R(s)$ denote the multi-set of all coordinates of $D(R)$ which do not exceed $s-1$ then
\begin{equation*}
\Big|\bigcup_{j=1}^t\ \{a_j+bi: i\in [0, s-1]\}\Big|=st-\sum_{R\in \mathcal{R}}\ \sum_{r\in M_R(s)}\ (s-r).
\end{equation*}
\end{prp}
We will refer to the sum
\[
\sum_{R\in \mathcal{R}}\ \sum_{r\in M_R(s)}\ (s-r)
\]
as the \emph{defect} of the set $\bigcup_{j=1}^t\ \{a_j+bi: i\in [0, s-1]\}$; we note that the defect is always non-negative, and it is equal to 0 if and only if the set $\{R\in \mathcal{R}: M_R(s)\not= \emptyset\}$ is empty. Returning to the $2m$-tuple $(\textbf{a}, \textbf{b})$ with which we began this section, we let $\Delta_i$ denote the defect of the set $S_i$ defined by (7.1) and so deduce the following proposition from Proposition 7.1:
\begin{prp}
\label{prp2}
If $\alpha$ denotes the sum of the cardinalities of the sets $S_1,\dots, S_k$ then
\begin{equation*}
\alpha=ms-\sum_{i=1}^k \Delta_i.
\]
 \end{prp} 
 
 \emph{Remark}. Suppose that $m=1$ in the definition of $AP(\textbf{a}, \textbf{b};s)$. In this case,  the technical details of the proof of Theorem 6.1 simplify to such an extent that we can prove that $AP(a, b; s)$ has the universal pattern property and the residue and non-residue support properties (see the beginning of section 2 for the definition of these properties). In fact, if $\varepsilon$ is a choice of signs for $[0, s-1]$, one can easily modify the proof of 
Theorem 6.1 to show that as $p\rightarrow +\infty$, the cardinality of the set
\[
\{A\in AP(a, b; s)\cap 2^{[1,  p-1]}: (A, \varepsilon, p)\ \textrm{is a residue pattern of}\  p\}
\]
is asymptotic to $(b\cdot2^s)^{-1}p$. With a bit more effort, one can also prove that as $p\rightarrow +\infty$, the cardinality of the set
\[
\{A\in AP(a, b; s)\cap 2^{[1,  p-1]}: A\ \textrm{is a residue (respectively, non-residue) support set of}\ p\}
\]  
is asymptotic to $(b\cdot 2^{1+b(s-1)})^{-1}p$. We leave the details to the interested reader.
\section{An interesting class of examples}
 Let $k \in [2, +\infty)$. We will say that a $2k$-tuple $(\textbf{a}, \textbf{b})$  is \emph{admissible} if it satisfies the following two conditions:
\begin{equation*}
\textrm{the coordinates of}\ \textbf{b}\ \textrm {are distinct, and}\tag{8.1},
\]
\begin{equation*}
a_ib_j-a_jb_i\not= 0\ \textrm{for}\ i\not= j\tag{8.2}.
\end{equation*}
If $s\in [1, +\infty)$ then it follows trivially from (8.1) that the parameter $\alpha$ in the statement of Theorem 6.1 for $AP(\textbf{a}, \textbf{b}; s)$ is $ks$, and so nothing more is needed for its calculation. However, when  $s\geq 2$ and $(\textbf{a}, \textbf{b})$ is admissible, we will indicate how overlap diagrams can be used to calculate the parameter $e$ in the statement of Theorem 6.1, in fact in a manner very similar to the way that they were used before in the proof of Proposition 7.1. We will then use this calculation of $e$ to illustrate precisely how Theorem 6.1 operates in some concrete situations.

Begin by letting $q_i=a_i/b_i$ for $i\in [1, k]$; without loss of generality, we suppose that the coordinates of $\textbf{a}$ and $\textbf{b}$ are indexed so that $q_i<q_{i+1}$ for each $i\in [1, k-1]$. Consider now the set $Q(\textbf{a}, \textbf{b})$ of all elements $(i, j)$ of $[1, k]\times [1, k]$ such that $i\not= j$ and $b_ib_j$ divides $a_ib_j-a_jb_i$, with  quotient $q(i, j)$, say. Condition (8.2) guarantees that if $(i, j)\in Q(\textbf{a}, \textbf{b})$ then $q(i, j)\not= 0$. Because $\Lambda(\mathcal{K})$ is empty if and only if $Q(\textbf{a}, \textbf{b})$ contains no elements $(i, j)$ such that $|q(i, j)|\leq s-1$, we need only calculate $e$ when  $Q(\textbf{a}, \textbf{b})$ contains elements of this type. Thus, suppose that this is so.

Let $\pi$ denote the canonical projection of $[1, k]\times [1, k]$ onto its left factor. If $(i, j)\in \pi(Q(\textbf{a}, \textbf{b}))\times \pi(Q(\textbf{a}, \textbf{b}))$ and we declare that $i\simeq j$ if either $i=j$ or $(i, j) \in Q(\textbf{a}, \textbf{b})$, then $\simeq$ defines an equivalence relation on $\pi(Q(\textbf{a}, \textbf{b}))$. 

We will now construct a series of overlap diagrams in a manner very similar to the procedure that we used in the proof of Proposition 7.1. Let $F$ be an equivalence class of the equivalence relation $\simeq $ such that $|q(i, j)|\leq s-1$ for some $(i, j)\in F\times F$. We note that the set $\{q_i: i\in F\}$ is linearly ordered by the equations $q_i-q_j=q(i, j)$ for $i, j\in F$ with $i\not= j$. Next, consider the nonempty and pairwise disjoint family of all subsets $S$ of $\{q_i: i\in F\}$ such that $|S|\geq 2$ and $S$ is maximal with respect to the property that the distance between consecutive elements of $S$ does not exceed $s-1$. We  index the positive quotients $q(i, j)$ which implement the linear ordering in $S$ as $(q_S(i): i\in [1, |S|-1])$, and then let $\mathcal{D}(S)$ denote the overlap diagram of this $(|S|-1)$-tuple. Because $q_S(i)\leq s-1$ for all $i\in [1, |S|-1]$, $\mathcal{D}(S)$ consists of a single block.  

Using a suitable positive integer $v$, we index all of the sets $S$ that arise from all of the equivalence classes  in the previous construction as $S_1,\dots, S_v$ and then define the \emph{quotient diagram} of $(\textbf{a}, \textbf{b})$ to be the $v$-tuple of overlap diagrams $(\mathcal{D}(S_n): n\in [1, v])$. One can then prove that
\[
e=\sum_{n=1}^v\ \ \sum_{i=1}^{|S_n|-1} (s-q_{S_n}(i)). 
\]

Let $\mathcal{Q}$ denote the set of all subsets $Q$ of $\{q_1,\dots,q_k\}$ such that $|Q|\geq 2$ and $Q$ is maximal relative to the property that $w-z$ is an integer for all $(w, z)\in Q\times Q$. After linearly ordering the elements of each $Q\in \mathcal{Q}$, we let $D(Q)$ denote the $(|Q|-1)$-tuple of positive integers whose coordinates are the distances between consecutive elements of $Q$. Then if $M_Q(s)$ denotes the multi-set formed by the coordinates of $D(Q)$ which do not exceed $s-1$, the above formula for $e$ can be rewritten as 
\begin{equation*}
e=\sum_{Q\in \mathcal{Q}}\ \sum_{q\in M_Q(s)} \ (s-q)\tag{8.3}.
\]
Thus, in the same spirit and by the same method of Proposition 7.1, $e$ can be calculated solely by means of information obtained directly and straightforwardly from the set $\{q_1,\dots,q_k\}$.   
 
The quotient diagram $\mathcal{D}$ of $(\textbf{a}, \textbf{b})$ can be used to calculate the set $\Lambda(\mathcal{K})$ determined by $(\textbf{a}, \textbf{b})$ and hence also the associated signature of an allowable prime. In order to do that, let $S_1,\dots,S_v$ be the subsets of $\{q_1,\dots, q_k\}$ that determine the sequence of overlap diagrams $\mathcal{D}(S_1),\dots,\mathcal{D}(S_v)$ which constitute $\mathcal{D}$, and then find the subset $J_n$ of $[1, k]$ such that $S_n=\{q_j: j\in J_n\}$.The overlap diagram $\mathcal{D}(S_n)$ consists of $|J_n|$ rows, with each row containing $s$ points. If $i\in [1, |J_n|]$ then there is a unique element $j$ of $J_n$ such that the $i$-th element of $S_n$ is $q_j$; we now take $l\in [1, s]$ and label the $l$-th point of row $i$ in $\mathcal{D}(S_n)$ as $(j, l-1)$. If $\mathcal{C}$ denotes the set of all columns of $\mathcal{D}$ then we identify a column $C\in \mathcal{C}$ with the subset of $[1, k]\times [0, s-1]$ defined by  
\[
\{(i, j)\in [1, k]\times [0, s-1]: (i, j)\ \textrm{is the label of a point in}\ C\}.
\]
One can then show that if $\theta$ denotes the projection of $[1, k]\times [0, s-1]$ onto $[1, k]$ then
\begin{equation*}
\Lambda(\mathcal{K})=\bigcup_{C\in \mathcal{C}}\ \theta(\mathcal{E}(C)).\tag{8.4}
\]
When this formula for $\Lambda(\mathcal{K})$ is combined with (8.3), it follows that all of the data required for an application of Theorem 6.1 can be easily read off directly from the quotient diagram of $(\textbf{a}, \textbf{b})$.
 
Let $v\in [1, +\infty)$ and for each $n\in [1, v]$, let $\mathcal{D}(n)$ be a fixed but arbitrary overlap diagram with $k_n$ rows, $k_n\geq 2$, and gap sequence $(d(i, n): i\in [1, k_n-1])$, with no gap exceeding $s-1$. Let $k_0=0, k= \sum_n k_n$. We will now exhibit infinitely many admissible $2k$-tuples $(\textbf{a}, \textbf{b})$ whose quotient diagram is $\mathcal{D}=(\mathcal{D}(n): n\in [1, v])$. This is done by taking the $(k-1)$-tuple $(d_1,\dots, d_{k-1})$ in the following lemma to be 
\[
d_i=\left\{\begin{array}{cc}d\Big(i-\sum_0^n k_j, n+1\Big),\ \textrm{if}\ i\in \Big[1+\sum_0^n k_j, -1+\sum_0^{n+1} k_j\Big],\ n\in [0, v-1],\\
s,\ \textrm{elsewhere,}\end{array}\right.
\]
and then letting $(\textbf{a}, \textbf{b})$ be any $2k$-tuple obtained from the construction in the lemma.

\begin{lem}
\label{lem1}
For $k\in [2, +\infty)$, let $(d_1,\dots, d_{k-1})$ be a $(k-1)$-tuple of positive integers. Define $k$-tuples $(a_1,\dots, a_k), (b_1,\dots, b_k)$ of positive integers inductively as follows: let $(a_1, b_1)$ be arbitrary, and if $i>1$ and $(a_i, b_i)$ has been defined, choose $t_i\in [2, +\infty)$ and set
\[
a_{i+1}=t_i(a_i+d_ib_i),\ \ b_{i+1}=t_ib_i.
\]
Then
\[
a_ib_j-a_jb_i=\Big(\sum_{r=j}^{i-1}\ d_r\Big)b_ib_j,\ \textrm{for all}\  i>j.
\]
\end{lem}

We can also find infinitely many admissible $2k$-tuples $(\textbf{a}, \textbf{b})$  with the given quotient diagram  and  such that the set $\Pi_-$ determined by  $(\textbf{a}, \textbf{b})$ is empty. To do this, simply choose the integer $b_1$ and all subsequent $t_i$'s used in the above construction from Lemma 8.1 to be squares. This shows that there are infinitely many admissible $2k$-tuples with a specified quotient diagram which satisfy the hypothesis of Theorem 6.1($ii$)($b$). On the other hand, if  $b_1$ and all the subsequent $t_i$'s are instead chosen to be distinct primes, it follows that the $2k$-tuples determined in this way all have quotient diagram $\mathcal{D}$ and each have $\Pi_-$ of infinite cardinality, and so there are infinitely many admissible $2k$-tuples with specified quotient diagram which satisfy the hypothesis of Theorem 6.1($ii$)($c$). We also note that if $m\in [1, +\infty)$ and $(\textbf{a}, \textbf{b})$ is a fixed $2m$-tuple, one can easily find infinitely many ordered pairs $(a, b)\in [1, +\infty)\times [1, +\infty)$ such that $b\not= b_i$ and $bb_i$ does not divide $a_ib-ab_i$ for $i\in [1, m]$. Hence there are infinitely many admissible $2k$-tuples which satisfy the hypothesis of Theorem 6.1($i$).

With this cornucopia of examples in hand, for $\varepsilon\in \{-1, 1\}$, we let $c_\varepsilon(p)$ denote the cardinality of the set
\begin{equation*}
\{A\in AP(\textbf{a}, \textbf{b}; s)\cap 2^{[1, p-1]}: \chi_p(a)=\varepsilon, \ \textrm{for all}\ a\in A\},
\end{equation*} 
where $(\textbf{a}, \textbf{b})$ is admissible. We will now use the quotient diagram of $(\textbf{a}, \textbf{b})$, formulae (8.3), (8.4), and Theorem 6.1 to study how $(\textbf{a}, \textbf{b})$ determines the asymptotic behavior of $c_{\varepsilon}(p)$ in specific situations. We will illustrate how things work when $k=2$ and 3, and for when "minimal" or "maximal" overlap is present in the quotient diagram of $(\textbf{a}, \textbf{b})$.     

When $k=2$, there is only at most a single overlap of rows in the quotient diagram $\mathcal{D}$  of $(\textbf{a}, \textbf{b})$, and if, e.g., $a_1b_2-a_2b_1=qb_1b_2$ with $0<q\leq s-1$, then the quotient diagram looks like
\begin{center}
\begin{tabular}{ccccccccccccc}
$\cdot$&$\cdot$&$\cdot$&$\cdot$\ &\ $\cdot$\ &\ $\cdot$\ &\ $\cdot$\ &\ $\cdot$\ &&&&&\\
$\leftarrow$&$q$&$\rightarrow$&$\cdot$\ &\ $\cdot$\ &\ $\cdot$\ &\ $\cdot$\ &\ $\cdot$\ &\ $\cdot$\ &\ $\cdot$\ &\ $\cdot$\  &\ , \\
\end{tabular}
\end{center}
\vspace{0.5cm}
\noindent where $\alpha=2s$ and, because of (8.3), $e=s-q$. Formula  (8.4) shows that the signature of $p$ is $\{\chi_p(b_1b_2)\}$, and so we conclude from Theorem 6.1 that when $b_1b_2$ is a square,
\[
c_\varepsilon(p)\sim (b\cdot 2^{s+q})^{-1}p,\ \textrm{as}\ p\rightarrow +\infty,
\]
and when $b_1b_2$ is not a square, $\Pi_+$ is the set of all allowable primes $p$ such that $\{b_1, b_2\}$ is either a set of quadratic residues of $p$ or a set of quadratic non-residues of $p$, $\Pi_-$ is the set of all allowable primes $p$ such that $\{b_1, b_2\}$ contains a quadratic residue of $p$ and a quadratic non-residue of $p$,  
\[
c_\varepsilon(p)=0,\ \textrm{for all}\ p\ \textrm{in}\ \Pi_-,
\]
and as $p\rightarrow +\infty$ inside $\Pi_+$,
\[
c_\varepsilon(p)\sim (b\cdot 2^{s+q})^{-1}p.
\]

 When $k=3$ there are exactly three types of overlap possible in the quotient diagram of $(\textbf{a}, \textbf{b})$, determined, e.g., when either

 $(i)$ exactly one,

 $(ii)$ exactly two, or

 $(iii)$ exactly three
 
 \noindent of $b_1b_2$, $b_2b_3$, and $b_1b_3$ divide, respectively, $a_2b_1-a_1b_2, a_3b_2-a_2b_3$, and $a_3b_1-a_1b_3$ with positive quotients not exceeding $s-1$ in modulus.

In case $(i)$, with $a_2b_1-a_1b_2=qb_1b_2$, say, the block in the quotient diagram of $(\textbf{a}, \textbf{b})$ is formed by a single overlap between rows 1 and 2, and this block looks exactly like the overlap diagram that was displayed for $k=2$ above. It follows that the conclusions from (8.3), (8.4), and Theorem 6.1 in case $(i)$ read exactly like the conclusions in the $k=2$ case described before, except that the exponent of the power of $1/ 2$ in the coefficient of $p$ in the asymptotic approximation is now $2s+q$ rather than $s+q$.   

In case $(ii)$, with $a_2b_1-a_1b_2=qb_1b_2$ and $a_3b_2-a_2b_3=rb_2b_3$, say, the block in the quotient diagram is formed by an overlap between rows 1 and 2 and an overlap between rows 2 and 3. Hence it looks like 
\vspace{0.5cm}
\begin{center}
\begin{tabular}{cccccccccccccccccc}
$\cdot$&$\cdot$&$\cdot$&$\cdot$\ &\ $\cdot$\ &\ $\cdot$\ &\ $\cdot$\ &\ $\cdot$&& &&&&&&&&\\
$\leftarrow$&$q$&$\rightarrow$&$\cdot$\ &\ $\cdot$\ &\ $\cdot$\ &\ $\cdot$\ &\ $\cdot$\ &\ $\cdot$\ &\ $\cdot$\ &\ $\cdot$\  &\  \\
&&&$\leftarrow$&&$r$&&$\rightarrow$\ &\ $\cdot$\ &\ $\cdot$\ &\ $\cdot$\ &\ $\cdot$\ \ &\  $\cdot$\ \ &\ $\cdot$\ \ &\ $\cdot$\ \ &\ $\cdot$\  &\ , \\
\end{tabular}
\end{center}
\vspace{0.5cm}
\noindent where $\alpha=3s$, and, because of (8.3) and (8.4), $e=2s-q-r$ and the signature of $p$ is $\{\chi_p(b_1b_2), \chi_p(b_2b_3)\}$. We hence conclude from Theorem 6.1 that if $b_1b_2$ and $b_2b_3$ are both squares then
\begin{equation*}
c_\varepsilon(p)\sim (b\cdot 2^{s+q+r})^{-1}p\ \textrm{as}\ p\rightarrow +\infty.\tag{8.5}
\end{equation*}  
On the other hand, if either $b_1b_2$ or $b_2b_3$ is not a square then $\Pi_+ $ consists of all allowable primes $p$ such that $\{b_1, b_2, b_3\}$  is either a set of quadratic residues of $p$ or a set of quadratic non-residues of $p$, $\Pi_-$ consists of all allowable primes $p$ such that  $\{b_1, b_2, b_3\}$ contains a quadratic residue of $p$ and a quadratic non-residue of $p$,
\begin{equation*}
c_\varepsilon(p)=0,\  \textrm{for all}\ p\in \Pi_-,\ \textrm{and}\tag{8.6}
\end{equation*} 
\begin{equation*}
c_\varepsilon(p)\sim (b\cdot 2^{s+q+r})^{-1}p\ \textrm{as}\ p\rightarrow +\infty\ \textrm{inside}\ \Pi_+.\tag{8.7}
\end{equation*} 

In case $(iii)$, with the quotients $q$ and $r$ determined as in case $(ii)$, and, in addition, $a_3b_1-a_1b_3=tb_1b_3$, say, the block in the quotient diagram is now formed by an overlap between each pair of rows, and so the diagram looks like
\vspace{0.5cm}
\begin{center}
\begin{tabular}{cccccccccccccccccc}
$\cdot$&$\cdot$&$\cdot$&$\cdot$\ &\ $\cdot$\ &\ $\cdot$\ &\ $\cdot$\ &\ $\cdot$&& &&&&&&&&\\
$\leftarrow$&$q$&$\rightarrow$&$\cdot$\ &\ $\cdot$\ &\ $\cdot$\ &\ $\cdot$\ &\ $\cdot$\ &\ $\cdot$\ &\ $\cdot$&$\cdot$&  \\
&&&$\leftarrow$&$r$&$\rightarrow$\ &\ $\cdot$\ &\ $\cdot$\ &\ $\cdot$\ &\ \ $\cdot$\ \ &\  $\cdot$\ \ &\ $\cdot$\ \ &\ $\cdot$\ \ &\ $\cdot$\  &\ , \\
\end{tabular}
\end{center}  
\vspace{0.5cm}

\noindent where $\alpha= 3s, e=2s-q-r$, and the signature of $p$ is $\{\chi_p(b_1b_2), \chi_p(b_1b_3), \chi_p(b_2b_3)\}$. In this case, the asymptotic approximation (8.5) holds whenever $b_1b_2, b_1b_3$, and $b_2b_3$ are all squares, and when at least one of these integers is not a square, $\Pi_+$ and $\Pi_-$ are determined by $\{b_1, b_2, b_3\}$ as before and (8.6) and (8.7) are valid.

\emph{Minimal overlap}. Here we take the quotient diagram to consist of a single block with gap sequence $(s-1, s-1,\dots,s-1)$, so that the overlap between rows is as small as possible: a typical quotient diagram for $k=5$ looks like
\vspace{0.5cm}
\begin{center}
\begin{tabular}{cccccccccccccccc}
$\cdot$&$\cdot$&$\cdot$&$\cdot$&&&&&&&&&&&&\\
&&&$\cdot$&$\cdot$&$\cdot$&$\cdot$&&&&&&&&&\\
&&&&&&$\cdot$&$\cdot$&$\cdot$&$\cdot$&&&&&&\\
&&&&&&&&&$\cdot$&$\cdot$&$\cdot$&$\cdot$&\\
&&&&&&&&&&&&$\cdot$&$\cdot$&$\cdot$&$\cdot$\ \   .\\
\end{tabular}
\end{center}
\vspace{0.5cm}
\noindent Here $\alpha=ks$, $e=k-1$, and the signature of $p$ is $\{\chi_p(b_ib_{i+1}): i\in [1, k-1]\}$. Hence via Theorem 6.1 , if $b_ib_{i+1}, i\in [1, k-1]$, are all squares then
\begin{equation*}
c_\varepsilon(p)\sim (b\cdot 2^{1+k(s-1)})^{-1}p\ \textrm{as}\ p\rightarrow +\infty,
\end{equation*}  
and if at least one of those products is not a square, then $\Pi_+ $ consists of all allowable primes $p$ such that $\{b_1, \dots, b_k\}$  is either a set of quadratic residues of $p$ or a set of quadratic non-residues of $p$, $\Pi_-$ consists of all allowable primes $p$ such that  $\{b_1, \dots, b_k\}$ contains a quadratic residue of $p$ and a quadratic non-residue of $p$,
\begin{equation*}
c_\varepsilon(p)=0,\  \textrm{for all}\ p\in \Pi_-,\ \textrm{and}\tag{8.8}
\end{equation*} 
\begin{equation*}
c_\varepsilon(p)\sim (b\cdot 2^{1+k(s-1)})^{-1}p\ \textrm{as}\ p\rightarrow +\infty\ \textrm{inside}\ \Pi_+.
\end{equation*} 

\emph{Maximal overlap} ($k\geq 3$). Here we take the quotient diagram to consist of a single block with gap sequence $(1,1,\dots,1)$, so that the overlap between each pair of rows is as large as possible: the diagrams for $k=3, 4,$ and 5 look like
\vspace{0.5cm}
\begin{center}
\begin{tabular}{ccccccccccccccccccccccc}
$\cdot$&$\cdot$&$\cdot$&&&&$\cdot$&$\cdot$&$\cdot$&$\cdot$&&&&&$\cdot$&$\cdot$&$\cdot$&$\cdot$&$\cdot$&&&&\\
&$\cdot$&$\cdot$&$\cdot$&&&&$\cdot$&$\cdot$&$\cdot$&$\cdot$&&&&&$\cdot$&$\cdot$&$\cdot$&$\cdot$&$\cdot$&&&\\
&&$\cdot$&$\cdot$&$\cdot$&&&&$\cdot$&$\cdot$&$\cdot$&$\cdot$&&&&&$\cdot$&$\cdot$&$\cdot$&$\cdot$&$\cdot$&&\\
&&&&&&&&&$\cdot$&$\cdot$&$\cdot$&$\cdot$&&&&&$\cdot$&$\cdot$&$\cdot$&$\cdot$&$\cdot$ \\
&&&&&&&&&&&&&&&&&&$\cdot$&$\cdot$&$\cdot$&$\cdot$&$\cdot$\ \  .\\
\end{tabular}
\end{center}
\vspace{0.5cm}
\noindent We have in this case that $k=s$, $\alpha=k^2$, $e=(k-1)^2$, and the signature of $p$ is
\[
\Big\{\chi_p\Big(\prod_{i\in I} b_i\Big): \emptyset\not= I\subseteq [1, k], |I|\ \textrm{even}\ \Big\}.
\]
 Hence if $\prod_{i\in I} b_i$ is a square for all nonempty subsets $I$ of $[1, k]$ of even cardinality then
\begin{equation*}
c_\varepsilon(p)\sim (b\cdot 2^{2k-1})^{-1}p\ \textrm{as}\ p\rightarrow +\infty,
\end{equation*}  
and if one of these products is not a square then $\Pi_+$ and $ \Pi_-$ are determined by  $\{b_1, \dots, b_k\}$ as before, (8.8) holds, and 
\begin{equation*}
c_\varepsilon(p)\sim (b\cdot 2^{2k-1})^{-1}p\ \textrm{as}\ p\rightarrow +\infty\ \textrm{inside}\ \Pi_+.
\end{equation*}

It follows from our discussion after the proof of Theorem 6.1 in section 6 that an increase in the number of overlaps between rows in $\mathcal{D}$ leads to an increase in the number of elements of $AP(\textbf{a}, \textbf{b}; s)\cap 2^{[1, p-1]}$ that are sets of quadratic residues or non-residues of $p$, and these examples now verify that principle quantitatively.  In order to see this explicitly, note first that Lemma 8.1 can be used to generate examples where $b$ always takes the same value. Hence we may assume in the discussion to follow that the value of $b$ is constant in each set of examples, and so the only parameter that is relevant when comparing asymptotic approximations to $c_{\varepsilon}(p)$ is the exponent of the power of $1/2$ in the coefficient of that approximation. When $k=3$ there are, respectively, 1, 2, and 3 overlaps between rows in cases ($i$), ($ii)$, and $(iii)$. It follows that $q+r\geq s$ in case $(ii)$ and $q+r<s$ in case $(iii)$. Hence the exponent in the power of $1/ 2$ that occurs in the asymptotic approximation to $c_\varepsilon(p)$ is greater than $2s$ in case $(i)$, is at least $2s$ in case ($ii)$, and is less than $2s$ in case ($iii)$. If we also take $k=s$ when there is minimal overlap in $\mathcal{D}$ and compare that to what happens when there is maximal overlap in $\mathcal{D}$, we see that the exponent in the power of  $1/ 2$ that occurs in the asymptotic approximation of $c_\varepsilon(p)$ is quadratic in $k$, i.e., $ k^2-k+1$, in the former case, but only linear in $k$, i.e., $ 2k-1$, in the latter case.

\section{Some problems worthy of further study}
A natural companion to the sequences in arithmetic progression are the sequences in {\it geometric} progression, hence one may inquire about the distribution of quadratic residues and non-residues among the geometric progressions in $[1,p-1]$. If $(a,b)\in [1,+\infty)\times [2,+\infty)$, then the only residue pattern $\varepsilon$ of a prime $p$ that the geometric progression $\{ab^n: n\in [0,+\infty)\}$ can have is $\varepsilon \equiv 1$ (when $\chi_p(a)=\chi_p(b)=1), \; \varepsilon \equiv -1$ (when $\chi_p(a)=-1$ and $ \chi_p(b)=1$) or a pattern of alternating signs (when $\chi_p(b)=-1$). In order to move beyond this rather trivial situation, we say that a sequence is in {\it almost geometric progression} if it has the form
\[
\{c+ab^n: \; n\in [0,+\infty )\}
\]
for a fixed element $(a,b,c)$ of $[1,+\infty )\times [2,+\infty )\times [1,+\infty )$. If $s\in[1,+\infty)$ then one can ask, for example, does the set
\[
\left\{\{c+ab^{(n+i)}: \; i\in [0,s-1]\}: \; n\in [0,+\infty )\right\}
\]
have the universal pattern property, the residue support property, and/or the non-residue support property, and if so, what are the relevant asymptotics? We will address these questions and related ones in a forthcoming paper [18]. One very interesting aspect of this problem, among others, is the estimation of restricted Weil-type  sums with an ``exponential argument,'' i.e., sums of the form
\[
\sum^{[\log_b p]}_{x=0}\chi_p\left(f(b^x)\right),
\]
where $f$ is a polynomial as in Lemma 2.2 with $f(0)\neq 0$.

 Let $\Pi$ denote the set of all odd primes. If $(a,b,s,p)\in [0,+\infty)\times [1,+\infty)\times[1,+\infty)\times \Pi$ and $\varepsilon$ is a choice of signs for $[1,s]$, we define
\[
n_0(p,s,\varepsilon )=\min \left\{n\in [1,+\infty ):\; \left(\left\{a+b(n+i):\; i\in [0,s-1]\right\}, \; \varepsilon,p\right)\text{ is a residue pattern of }p\right\},
\]
\[
\ \ n_1^+(p,s) =\min \left\{n\in [1,+\infty ):\; \left\{a+b(n+i):\; i\in [0,s-1]\right\} \; \text{ is a residue support set  of }p\right\},
\]
\[
s_0^+(p)=\max \left\{n\in [1,+\infty ): \; \text{ there exists } S\in AP(a,b;\, n)\text{ such that }S\subseteq R(p)\right\},
\]
\[
s_0^-(p)=\max \left\{n\in [1,+\infty ): \; \text{ there exists } S\in AP(a,b;\, n)\text{ such that }S\subseteq NR(p)\right\},
\]
\[
s^+_1(p)=\max \{n \in [1, +\infty ):\text{ there exists }S\in AP(a,b;\, n)\text{ such 
that }S\text{ is a residue support set of }p\},
\]
\[
q^{\pm}_0(s)=\min \left\{ q\in \Pi :\: s^{\pm}_0(q)=s\right\},
\]
\begin{eqnarray*}
Q^{\pm}_0(s)&=&\max \left\{q\in\Pi:\; s^{\pm}_0(q)=s\right\},\\
q_1^+(s)&=&\min \left\{ q\in \Pi :\: s_1^+(q)=s\right\},\\
Q_1^+(s)&=&\max \left\{q\in\Pi:\; s_1^+(q)=s\right\}.
\end{eqnarray*}
We also define $n^-_1, s_1^-,q^-_1,$ and $Q^-_1$ in a similar fashion, using non-residue support sets in place of residue support sets.

It is a consequence of the remark  at the end of section 7 that each of these functions is finite-valued for all $p$ sufficiently large, for all $s\in [1,+\infty )$ and for all choice of signs $\varepsilon$. The function $n_0$ locates the first occurrence of a residue pattern of specified length and type of a fixed prime in consecutive terms of $AP(a,b)$, and $s^+_0$ (respectively, $s^-_0$) is the longest length of a set of quadratic residues (respectively, non-residues) of a fixed prime located in consecutive terms of $AP(a,b)$. The parameter $q^+_0 (s)$ (respectively, $q^-_0 (s)$)  is the smallest prime which has a set of quadratic residues (respectively, non-residues) of length $s$  in consecutive terms of $AP(a,b)$ and $Q^+_0(s)$ (respectively, $Q^-_0(s)$) is the largest prime that does not have a set of quadratic residues (respectively, non-residues) of length $s+1$ in consecutive terms of $AP(a,b)$. Similar descriptions of the functions $n^{\pm}_1, s_1^{\pm},q^{\pm}_1,$ and $Q^{\pm}_1$ also hold. We note finally that all of these functions have analogs that are defined by replacing elements in $AP(a,b;s)$ and $AP(a,b;n)$ by elements of $AP(b;s)$ and $AP(b;n)$ in the above definitions wherever appropriate, with descriptions similar to those just given also valid for these analogs. 

Our last problem calls for a detailed study of the functions defined above, with an emphasis on nontrivial estimates and asymptotics in terms of the relevant variables. Such a study would uncover much interesting information about the fine structure of quadratic residues and non-residues in arithmetic progression. For  $AP(0,1)=[0,+\infty )$, Burgess [4] obtained a good estimate of $s_0^+(p)$ and $s^-_0(p)$. More recently, Buell and Hudson [1] and Hudson [10] have obtained noteworthy results on the behavior of $n_0,s^{\pm}_0,q^{\pm}_0$, and $Q^{\pm}_0$ for $AP(0,1)$ in certain special cases. What transpires for more general arithmetic progressions is completely open, as far as we know.

\textsc{acknowledgement}. I am deeply grateful to my dear wife Linda, whose encouragement and wise council during some discouraging moments kept me focused effectively on the task at hand.

\end{document}